\documentclass[11pt, a4paper, reqno]{amsart}

\usepackage[utf8]{inputenc}
\usepackage[left=30mm, right=30mm, top=25mm, bottom=25mm]{geometry}
\usepackage{float}
\usepackage{dsfont}
\usepackage{enumitem}
\usepackage{mathrsfs}
\usepackage{amsmath, amssymb, amsthm, bm}
\usepackage[colorlinks=true, linkcolor=blue, filecolor=blue, citecolor=red, urlcolor=blue, colorlinks=true]{hyperref}
\makeatletter
\newcommand{\labelleditem}[2]{%
  \item[#1]%
  \phantomsection
  \def\@currentlabel{#1}%
  \label{#2}%
}
\makeatother
\usepackage[capitalize,noabbrev]{cleveref}
\usepackage{tikz}
\usepackage{tikz-cd}
\usepackage{pgfplots}
\usetikzlibrary{patterns}
\usetikzlibrary{intersections, pgfplots.fillbetween}
\pgfplotsset{compat=1.18} 
\pgfdeclarelayer{pre main}

\newtheorem{theorem}{Theorem}[section]
\newtheorem{proposition}[theorem]{Proposition}
\newtheorem{lemma}[theorem]{Lemma} 
\newtheorem{corollary}[theorem]{Corollary}

\theoremstyle{definition}
\newtheorem{definition}[theorem]{Definition}
 
\newtheorem{example}[theorem]{Example}
\newtheorem{remark}[theorem]{Remark}

\newcommand{\rank}{\mathrm{rank}}

\newcommand{\supp}{\mathrm{supp}}
\newcommand{\trop}{\mathrm{trop}}

\newcommand{\val}{\mathrm{val}}
\newcommand{\Mat}{\mathrm{Mat}}
\newcommand{\Kcal}{\mathcal{K}}
\newcommand{\Ocal}{\mathcal{O}}
\newcommand{\Log}{\mathrm{Log}}
\newcommand{\scrA}{\mathscr{A}}
\newcommand{\GL}{\mathrm{GL}}
\newcommand{\PP}{\mathbb{P}}
\newcommand{\T}{\mathbb{T}}
\newcommand{\Q}{\mathbb{Q}}
\newcommand{\R}{\mathbb{R}}
\newcommand{\C}{\mathbb{C}}

\title{How to tropicalize a non-Archimedean lattice}
\author{Yassine {Elmaazouz}}
\address{Yassine {Elmaazouz} (Caltech)}
\email{maazouz@caltech.edu}
\date{\today}


\begin{document}
    
    \begin{abstract}
        We present a method for computing the tropicalization of a lattice
        over the valuation ring of a non-Archimedean field. 
        Our construction relies on the \emph{entropy polynomial}, a tropical
        polynomial whose coefficients are determined by the valuations of
        minors of a matrix representing the lattice. 
        We show that the tropicalization of the lattice is the support of a
        polyhedral complex defined through this polynomial. 
        Moreover, we prove that this tropicalization is a finitely generated
        tropical semimodule over the nonnegative tropical scalars and provide
        explicit generators in terms of the entropy vector of the lattice. 
        Finally, we extend our construction to valuated bimatroids,
        establishing analogous results beyond the realizable setting.
    \end{abstract}

\maketitle

\section{Introduction}

    Let $K$ be the field of \emph{generalized Puiseux series} with complex coefficients. This is the field of formal series in one variable $t$ of the form
    \[
        \bm{f} = \sum_{\alpha \in \R} a_\alpha \ t^{\alpha }, \quad \text{where } \ \ a_\alpha \in \C \quad \text{for all } \alpha \in \R,
    \]
    such that the \emph{support} of $\bm{f}$,
    \(
        \supp(\bm{f}) := \{ \alpha \in \R \colon a_\alpha \neq 0 \} \subset \R,
    \) 
    is either finite or its only accumulation point is $+ \infty$. This is a non-archimedean valued field with valuation
    \[
        \val: K \to \overline{\R}, \quad \bm{f} \mapsto  \min\big( \supp(\bm{f})\big),
    \]
    where $\overline{\R} := \R \cup \{\infty\}$. Here the minimum of the empty set is $\infty$ by convention, so $\val(0) = \infty$. We note that the field $K$ is algebraically closed, see \cite[Theorem 6]{Markwig10}. We denote by $\Ocal_{K}$ its valuation ring 
    \[
        \Ocal_{K} := \Big\{ \bm{f} \in K  \colon  \val(\bm{f}) \geq 0 \Big\}.
    \]
    For tropical arithmetic on $\overline{\R}$, we use the notation 
    \[
        x \oplus y := \min(x,y) \quad \text{and} \quad x \odot y := x + y \qquad \text{whenever } x,y \in \overline{\R}.
    \] 
    The triplet $(\overline{\R}, \oplus, \odot)$ is the tropical semiring with the min-plus algebra. 
    We define tropical addition and scalar multiplication for vectors coordinatewise by setting $x \oplus y := (x_i \oplus y_i)_{i=1}^{n}$ and $\lambda \odot x := (\lambda \odot x_i)_{i=1}^{n}$ whenever $\lambda \in \overline{\R}$ and $x,y \in \overline{\R}^n$. 
    Finally we write $\overline{\R}_{\geq 0} := \R_{\geq 0} \cup \{\infty\}$.

    \subsection{Non-archimedean lattices} 
    
    A \emph{lattice} $L \subset K^n$ is a finitely generated submodule of $K^n$ over the valuation ring $\Ocal_{K}$ i.e., there exist $A_1, \dots, A_r \in K^n$ such that
    \[
        L = \Ocal_{K} A_1 + \dots + \Ocal_{K} A_r.
    \]
    We recall that lattices of rank $1 \leq r \leq n$ are in one-to-one correspondence with left orbits in $  \GL_r(\Ocal_{K}) \backslash \Mat^{\circ}_{r \times n}(K)$, where $\Mat^{\circ}_{r \times n}(K)$ is the set of $r \times n$ matrices of full rank. That is, any rank $r$ lattice $L \subset K^n$ can be written as
    \[
        L = \Ocal_{K}^{r} A  := {\rm rowspan}_{\Ocal_{K}}(A), \quad \text{where } A \in \Mat_{r \times n}^{\circ}(K)
    \]
    and the matrix $A$ is unique up to left multiplication by an element of $\GL_r(\Ocal_{K})$. To any lattice $L \subset K^n$, we can associate a linear subspace $V := L \otimes_{\Ocal_{K}} K$ of $K^n$. Throughout this article, without loss of generality, we make the assumption that the vector space $V$ is not contained in any of the coordinate hyperplanes of $K^n$. If not, we can simply replace $K^n$ with a smaller linear subspace.
    
    While the tropicalization of linear spaces is well established and widely studied, see \cite[Chapter 4]{MS2015} and \cite{Speyer,Rincon,Hampe}, the tropicalization of lattices, to the author’s knowledge, is less understood. In this paper, we initiate the study of the \emph{tropicalization} $\trop(L)$ of a lattice $L \subset K^n$ which we define as follows.
    
    \begin{definition}
    Let $L \subset K^{n}$ be a lattice of rank $r$. We define its tropicalization by
    \begin{equation}\label{eq:defOfTropL}
        \trop(L)  := \Big\{ \big(\val(\bm{f}_1), \dots, \val(\bm{f}_n)\big) : \ (\bm{f}_1, \dots, \bm{f}_n) \in L  \Big\} \subset \overline{\R}^n.
    \end{equation}  
    See Figures~\ref{fig:tropL1} and \ref{fig:tropL2} for examples.
    \end{definition}

    \begin{remark}
        The tropicalization defined in \eqref{eq:defOfTropL} is taken with respect to the standard basis $e_1, \dots, e_n$ of $K^n$. In general, if $B \in \GL_n(K)$ is a matrix whose rows $b_1, \dots, b_n$ form an ordered basis of $K^n$, we can define a tropicalization $\trop_B$ as follows:
        \[
            \trop_B(L) := \Big\{ \big(\val(\bm{f}_1), \dots, \val(\bm{f}_n)\big) \colon \bm{f}_1 b_1  + \dots + \bm{f}_n b_n \in L  \Big\},
        \]
        and for any lattice $L$ we have $\trop_B(L) = \trop(L \ B^{-1})$. So it suffices to compute $\trop(L)$ for any lattice $L$ in the standard basis $e_1, \dots, e_n$ of $K^n$.
    \end{remark}

    \subsection{The entropy vector of a lattice}
    
    We briefly recall the definition of the \emph{entropy map} from \cite{EL22}. Let $L$ be a lattice in $K^n$ of rank $1 \leq r \leq n$, spanned over $\Ocal_{K}$ by the rows of a full-rank $r \times n$ matrix $A$. The \emph{entropy vector} of $L$ is the vector $h(L) = (h_J(L))$ with entries in $\overline{\R}$ indexed by subsets $J$ of $[n] := \{1,2,\dots,n\}$. Its entries are given by:
    \[
        h_J(L) := \min_{\substack{I \subset [r] \\  |I| = |J|}} \val(\det(A_{I \times J})).
    \]
    The minimum is taken with respect to all subsets $I$ of $[r]$ such that $I$ and $J$ have the same size. If no such $I$ exists, the value of $h_J(L)$ is set to $\infty$, and if $J = \varnothing$ we set $h_{\varnothing}(L) = 0$. We note that the entropy vector $h(L)$ encodes the tropical Pl\"ucker vector of the linear space $\PP(V) \subset \PP^{n-1}$ where $V := L \otimes_{\Ocal_K} K$ since
    \[
        h_J(L) = \val(\det(A_{[r] \times J})) \quad \text{for all } J \in \binom{[n]}{r}.
    \]
    
    \begin{definition}
        The \emph{entropy polynomial} of a lattice $L \subset K^n$ is the tropical polynomial defined as follows:
        \begin{equation}\label{eq:tropPoly}
            \varphi_L(v) := \max_{J \subset [n]} \left( v_J  - h_J(L)   \right),
        \end{equation}
        where $v_J := \sum_{j \in J} v_j$ for any $\varnothing \neq J \subset [n]$ and by convention $v_{\varnothing} = 0$.
    \end{definition}

    \begin{example}
        Let $L$ be the lattice spanned over $\Ocal_{K}$ by the rows of the matrix
        \[
            \begin{bmatrix}
                1 & 1 & 1 \\ 0 & t & t^2
            \end{bmatrix}.
        \]
        The entropy vector of $L$ is 
        \[
            h_{\varnothing} = 0, \quad  h_1 = h_2 = h_3 = 0, \quad h_{12} = h_{23} = 1, \quad h_{13} = 2, \quad \text{and} \quad h_{123} = \infty.
        \]
    \end{example}

    \begin{remark} It was shown in \cite{EL22} that the entropy vector $h(L) = (h_J(L))$ does not depend on the choice of the matrix $A$. Moreover, using the Cauchy-Binet formula, we can see that there exists a generic matrix $U \in \GL_r(\Ocal_{K})$ such that 
        \[
            h_J(L) = \frac{1}{2} \val(\det(S_{J \times J})), \qquad \text{where } S := (UA)^\top UA.
        \]
        See \cite[Lemma 3.4]{ARVY}. So the set of entropy vectors
        \[
            \mathcal{H}_{n} := \Big\{ h(L) \colon L \subset K^n \text{ is a lattice} \Big\},
        \]
        can be thought of, up to a factor of two, as arising from the valuations of principal minors of symmetric matrices. 
        We also note that for any lattice $L \subset K^n$, the entropy vector $h(L) = (h_J(L))_{J \subset [n]}$ is supermodular, i.e.,
        \[
            h_I(L) + h_J(L) \leq h_{I \cap J}(L) + h_{I \cup J}(L), \quad \text{for any } I, J \subset [n],
        \]
        see \cite[Theorem 4.2]{EL22}.
    \end{remark}

    \begin{remark}
        The terminology \emph{entropy vector} comes from
        the connection with multivariate Gaussian distributions
        over non-Archimedean local fields, studied in
        \cite{EL22}. Indeed, suppose that ${\mathbb K}$ is a local
        field and $\Ocal_{\mathbb K}$ its valuation ring with residue field of cardinality $q$, and let $\mu_J$ denote the Haar measure on ${\mathbb K}^J$
        normalized by $\mu_J(\Ocal_{\mathbb K}^J)=1$.
        If $\pi_J:\mathbb{K}^n \to \mathbb{K}^J$ denotes the coordinate
        projection, then
        \[
            h_J(L) = -\log_q\big(\mu_J(\pi_J(L))\big).
        \]
        Thus, whenever $h_J(L)<\infty$, the quantity
        $-h_J(L)$ is the differential entropy, with
        logarithms to base $q$, of the uniform probability
        measure on $\pi_J(L)$.
        The definition of $h_J(L)$ in terms of valuations
        of minors extends this construction to arbitrary
        non-Archimedean valued fields, even when the field does 
        not admit a Haar measure, which is the case for generalized Puiseux series.
    \end{remark}

    \begin{definition}\label{def:SigmaL}
    The tropical polynomial $\varphi_L$ subdivides $\R^n$ into polyhedra which are its regions of linearity. We denote by $\Kcal_L$ the $n$-dimensional polyhedral complex resulting from this subdivision.  We inductively define a labeling $\ell \colon \Kcal_L \to 2^{[n]}$ on the cells of $\Kcal_L$ as follows:
    
    \begin{enumerate}[wide=40pt, leftmargin = 56pt]
        \item each maximal cell $\sigma$ of dimension $n$ corresponds to a monomial $v_J - h_J(L)$ for some $J \subset [n]$ and we set $\ell(\sigma) := J$.
        \item for any cell $\tau$ in $\Kcal_L$ we set: $\ell(\tau) := \bigcup\limits_{\substack{\sigma \in \Kcal_L \\ \tau \preceq \sigma}} \ell(\sigma)$.
    \end{enumerate}    
    
        The labeling $\ell$ defines a polyhedral complex $\Sigma_L \subset \Kcal_L$ as follows:
        \[
        \Sigma_L = \Big\{ \sigma \in \Kcal_L \colon  \ell(\sigma) = [n]\Big\}.
        \]
        From the definition of the labeling $\ell \colon \Kcal_L \to 2^{[n]}$ we readily see that
        if $\tau$ is a face of a cell $\sigma \in \Sigma_L$ then $\tau \in \Sigma_L$. 
        So $\Sigma_L$ is indeed a polyhedral complex. 
        Finally, we define the extended polyhedral complex $\overline{\Sigma}_L \subset \overline{\R}^{n}$ as the closure of $\Sigma_L$ in the tropical projective space
        \[
            \mathbb{TP}^{n} := \Big(\overline{\R}^n \setminus \{(\infty, \dots, \infty)\}\Big) \Big/ \ \overline{\R} \bm{1}_{n+1}.
        \]
        Here $\bm{1}_{n+1}$ is the all ones vector and $\overline{\R}^{n}$ is viewed as an affine chart of $\mathbb{TP}^{n}$ via the embedding
        \[
            \overline{\R}^{n} \to \mathbb{TP}^n, \quad (x_1, \dots, x_n) \mapsto [0:x_1:\dots:x_n].
        \]
    \end{definition}

  \begin{example}
        Let $L$ be the lattice in $K^2$ spanned by the rows of the matrix:
        \[
            A = \begin{bmatrix}
                    1-t^5 &   t+ t^3 \\ 3+t^2 & 3t+t^3
                \end{bmatrix}.
        \]
        The entropy vector of the lattice $L$ is:
        \[
         h_{\varnothing}(L) = 0, \quad h_1(L) = 0, \quad h_2(L) = 1, \quad \text{and} \quad \quad h_{12}(L) = 3.
        \]
        The polyhedral complex $\Kcal_L \subset \R^2$ induced by the polynomial $\varphi_L$ together with the labeling of its cells is depicted inFigure~\ref{fig:2}.
            \begin{figure}[H]
                \centering
                \scalebox{0.8}{
                \begin{tikzpicture}                    
                    \draw[blue, line width=2] (0,0)--(3,0);
                    \draw[blue, line width=2] (0,0)--(0,3);
                    \draw[blue, line width=2] (0,0)--(-2,-2);

                    \draw[blue, dashed, line width=2] (3,0)--(3,3);
                    \draw[blue, dashed, line width=2] (0,3)--(3,3);
                    
                    \draw (-2.5,-2.3) node {\tiny $\{1,2\}$};

                    \draw[-] (-5,-2)--(-2,-2);
                    \draw[-] (-2,-5)--(-2,-2);

                    \draw (-1.4, -2.1) node {$(0,1)$};
                    \draw (-0.5, 0.3) node {\tiny $\{1,2\}$};

                    \draw [fill=black]  (0,0)  circle (3pt);
                    \draw [fill=black] (-2,-2) circle (3pt);

                    \draw (2,2) node {$\{1,2\}$};
                    \draw (1.5,-3) node {$\{1\}$};
                    \draw (-4,1.5) node {$\{2\}$};
                    \draw (-4,-4) node {$\varnothing$};
                    
                    \draw (-0.5,3.2) node {\footnotesize $\{1,2\}$};
                    \draw (3.5,-0.2) node {\footnotesize $\{1,2\}$};
                    \draw (-5.5,-2.2) node {\footnotesize $\{2\}$};
                    \draw (-2.25,-5.5) node {\footnotesize $\{1\}$};

                    \draw (0.4,-0.4) node {$(2,3)$};
                    \draw (-1.2,-0.7) node {\footnotesize $\{1,2\}$};

                    \path[pattern color=blue, line width = 2, pattern=north west lines, opacity = 0.5] (3,0)-- (3,3) -- (0,3) -- plot[domain=1:3.5, smooth] (0,0) -- cycle;
                \end{tikzpicture}}
                \caption{The polyhedral complex $\Kcal_L$ (consisting of $2 + 5 + 4 = 11$ cells) 
                and its labeling $\ell \colon \Kcal_L \to 2^{[2]}$. The blue cells are the
                cells of $\Kcal_L$ with label $\{1,2\}$ i.e. the cells of the polyhedral 
                complex $\Sigma_L$. The dashed lines indicate where the unbounded cells of
                $\Sigma_L$ have been truncated for visualization.}
                \label{fig:2}
            \end{figure}
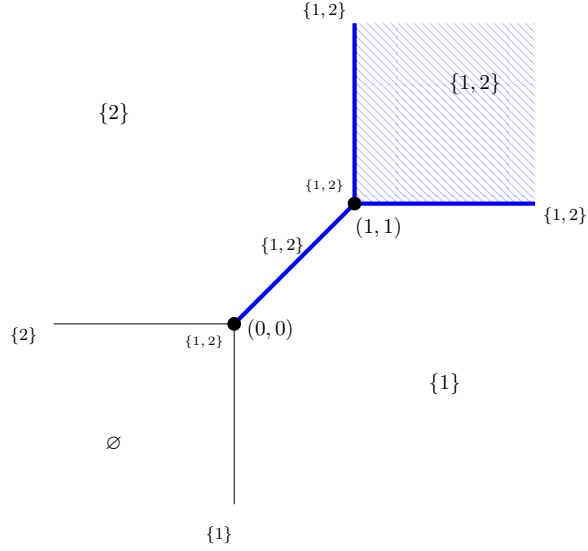
    \end{example}

\begin{remark}
    Let $L \subset K^n$ be a lattice of rank $r$, and for each $1 \leq k \leq r$, let
    \[
        h^{(k)}(L) :=
        \big(h_J(L)\big)_{J \in \binom{[n]}{k}}.
    \]
    Then the sequence
    \[
        \big(h^{(1)}(L),\dots,h^{(r)}(L)\big)
    \]
    defines a realizable valuated flag matroid whose
    underlying matroids are the successive truncations
    of its rank-$r$ matroid.

    Conversely, every realizable valuated flag matroid
    $(\mu^{(1)}, \dots, \mu^{(r)})$ whose underlying matroids are the 
    successive truncations of its rank-$r$ matroid arises in this way, up to
    an additive constant at each rank. More precisely, there exists a lattice 
    $L \subset K^n$ of rank $r$ and constants $c_1,\dots,c_r \in \R$ such that
    \[
       h_J(L) = \mu^{(k)}_J + c_k \qquad \text{for all } J \in \binom{[n]}{k}.
    \]
\end{remark}

    \subsection{Main results}
    We are now ready to state the main results of this paper.

    \begin{theorem}\label{thm:A}
        Let $L$ be a lattice of rank $r$ in $K^n$ that is not contained in any coordinate hyperplane of $K^n$. Then 
            \[
                \trop(L) = |\overline{\Sigma}_L| .
            \]
            The image of $|\Sigma_L|$ under the projection map $\R^n \to \T^{n-1} = \R^n/\R \bm{1}_n$ is equal to the tropical linear space $\trop(\PP(V))$ where $V := L \otimes_{\Ocal_K} K \subset K^n$. Consequently $\dim(\Sigma_L) = \rank(L)$.
    \end{theorem}
    
    \begin{remark} 
    The result in Theorem~\ref{thm:A} continues to hold true if one replaces $K$ with the field of Puiseux series $\C\{\!\{t\}\!\} := \bigcup_{N \geq 1} \C(\!(t^{1/N})\!)$ or the algebraic closure $\overline{\Q}_p$ of the field $\Q_p$ of $p$-adic numbers, except that we have to take the closure of $\trop(L)$ in the product topology on $\overline{\R}^{n}$.
    \end{remark}
    
    \begin{theorem}\label{thm:B}
        The tropical lattice $|\overline{\Sigma}_L|$ is a tropical semimodule over $\overline{\R}_{\geq 0}$. More precisely, for any $u,v \in |\overline{\Sigma}_L|$ and $\lambda\in\overline{\R}_{\geq 0}$ we have
        \[
            \lambda \odot u \in |\overline{\Sigma}_L| \quad \text{and} \quad u \oplus v \in |\overline{\Sigma}_L|.
        \]
        Moreover $|\overline{\Sigma}_L|$ is generated over $\overline{\R}_{\geq 0}$ by the set of vectors $\{ u_{J} \colon J \subsetneq [n], \ h_{J} \neq \infty \}$ in $\overline{\R}^n$ defined as follows 
        \[
            u_{J,j} := \begin{cases}
                        h_{Jj} - h_{J} & \text{if } j \not \in J \\ 
                             +  \infty & \text{otherwise} 
                      \end{cases}, \quad \text{where } Jj := J \cup \{j\}.
        \]
        More precisely we have
        \[
            |\overline{\Sigma}_L| = \bigoplus_{J \subsetneq [n], \ h_{J} \neq \infty} \overline{\R}_{\geq 0} \odot u_{J}.
        \]
    \end{theorem}
    
    \begin{remark}
         While it is tempting to say that $|\overline\Sigma_L|$ is the tropical cone generated by the vectors $u_I$, this is \emph{not} correct. In the literature, tropical cones are required to be closed under tropical addition and scalar multiplication with scalars in all of $\R$, see \cite{Joswig24}. Note that $|\Sigma_L|$ is closed under tropical addition and tropical scalar multiplication with scalars in $\R_{\geq 0}$.
    \end{remark}
    
    The following are a couple of examples of lattices and their tropicalizations.
    
    \begin{example}\label{ex:stdExample}
        Let $L$ be the lattice in $K^3$ spanned over $\Ocal_{K}$ by the rows of the following matrix:
        \begin{equation} \label{eq:matrixForExamlpe1}
            A = \begin{bmatrix}
                  1      &  1      &  1 \\
                  0      &  t^2    &  t \\
                  0      &  0      & t^2
            \end{bmatrix}.
        \end{equation}    
        The finite part $\trop(L)\cap\R^3$ is the support of the polyhedral
        complex $\Sigma_L$ depicted in Figure~\ref{fig:tropL1}.
        The entropy vector of $L$ is given by
     \[
        \begin{matrix}
                                  & h_{\varnothing}(L) = 0,   &\\
                     h_{1}(L) = 0, & h_{2}(L) = 0,           & h_{3}(L) = 0, \\
                    h_{12}(L) = 2, & h_{13}(L) = 1,          & h_{23}(L) = 1,\\
                                  & h_{123}(L) = 4.         &
        \end{matrix}
    \]
        \begin{figure}[ht]
        \centering
        \begin{minipage}{0.45\textwidth}
            \includegraphics[scale=0.26]{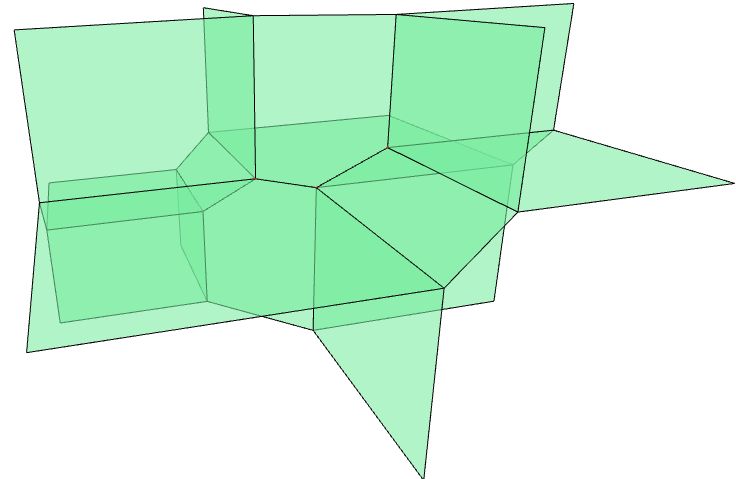}
        \end{minipage}
        \begin{minipage}{0.45\textwidth}
        \scalebox{0.4}{                
            \tikzset{every picture/.style={line width=0.75pt}}
            \begin{tikzpicture}[x=0.75pt,y=0.75pt,yscale=-1,xscale=1]
            \draw [line width=1.5]  [dash pattern={on 5.63pt off 4.5pt}]  (441.88,231.44) -- (496.73,191.14) ;
            \draw [line width=1.5]  [dash pattern={on 5.63pt off 4.5pt}]  (496.73,18.6) -- (496.73,191.14) ;
            \draw [line width=1.5]    (79.43,347.08) -- (134.28,298.73) ;
            \draw  [color={rgb, 255:red, 0; green, 0; blue, 0 }  ,draw opacity=0 ][fill={rgb, 255:red, 38; green, 119; blue, 249 }  ,fill opacity=0.33 ][line width=1.5]  (135.4,87.6) -- (276.89,105.31) -- (276.89,274.55) -- (222.04,314.85) -- (134.28,298.73) -- cycle ;
            \draw  [fill={rgb, 255:red, 0; green, 0; blue, 0 }  ,fill opacity=1 ][line width=1.5]  (74.93,347.08) .. controls (74.93,345.26) and (76.95,343.78) .. (79.43,343.78) .. controls (81.92,343.78) and (83.93,345.26) .. (83.93,347.08) .. controls (83.93,348.91) and (81.92,350.39) .. (79.43,350.39) .. controls (76.95,350.39) and (74.93,348.91) .. (74.93,347.08) -- cycle ;
            \draw  [fill={rgb, 255:red, 0; green, 0; blue, 0 }  ,fill opacity=1 ][line width=1.5]  (129.78,298.73) .. controls (129.78,296.9) and (131.8,295.42) .. (134.28,295.42) .. controls (136.77,295.42) and (138.78,296.9) .. (138.78,298.73) .. controls (138.78,300.55) and (136.77,302.03) .. (134.28,302.03) .. controls (131.8,302.03) and (129.78,300.55) .. (129.78,298.73) -- cycle ;
            \draw  [fill={rgb, 255:red, 0; green, 0; blue, 0 }  ,fill opacity=1 ][line width=1.5]  (217.54,314.85) .. controls (217.54,313.02) and (219.56,311.54) .. (222.04,311.54) .. controls (224.53,311.54) and (226.54,313.02) .. (226.54,314.85) .. controls (226.54,316.67) and (224.53,318.15) .. (222.04,318.15) .. controls (219.56,318.15) and (217.54,316.67) .. (217.54,314.85) -- cycle ;
            \draw [line width=1.5]    (134.28,298.73) -- (135.38,90.6) ;
            \draw [shift={(135.4,87.6)}, rotate = 90.3] [color={rgb, 255:red, 0; green, 0; blue, 0 }  ][line width=1.5]    (14.21,-4.28) .. controls (9.04,-1.82) and (4.3,-0.39) .. (0,0) .. controls (4.3,0.39) and (9.04,1.82) .. (14.21,4.28)   ;
            \draw  [fill={rgb, 255:red, 0; green, 0; blue, 0 }  ,fill opacity=1 ][line width=1.5]  (272.39,274.55) .. controls (272.39,272.73) and (274.41,271.25) .. (276.89,271.25) .. controls (279.38,271.25) and (281.39,272.73) .. (281.39,274.55) .. controls (281.39,276.38) and (279.38,277.86) .. (276.89,277.86) .. controls (274.41,277.86) and (272.39,276.38) .. (272.39,274.55) -- cycle ;
            \draw  [color={rgb, 255:red, 0; green, 0; blue, 0 }  ,draw opacity=0 ][fill={rgb, 255:red, 38; green, 119; blue, 249 }  ,fill opacity=0.33 ][line width=1.5]  (276.89,274.55) -- (496.73,191.14) -- (441.88,231.44) -- (222.04,314.85) -- cycle ;
            \draw  [color={rgb, 255:red, 0; green, 0; blue, 0 }  ,draw opacity=0 ][fill={rgb, 255:red, 38; green, 119; blue, 249 }  ,fill opacity=0.33 ][line width=1.5]  (276.89,274.55) -- (397.56,419.61) -- (342.71,463.21) -- (222.04,314.85) -- cycle ;
            \draw  [color={rgb, 255:red, 0; green, 0; blue, 0 }  ,draw opacity=0 ][fill={rgb, 255:red, 255; green, 17; blue, 0 }  ,fill opacity=0.12 ][line width=1.5]  (496.73,18.6) -- (617.4,107.25) -- (617.4,185.02) -- (617.4,336.2) -- (496.73,191.14) -- cycle ;
            \draw  [color={rgb, 255:red, 0; green, 0; blue, 0 }  ,draw opacity=0 ][fill={rgb, 255:red, 38; green, 249; blue, 50 }  ,fill opacity=0.25 ][line width=1.5]  (496.73,191.14) -- (617.4,336.2) -- (397.56,419.61) -- (276.89,274.55) -- cycle ;
            \draw  [color={rgb, 255:red, 0; green, 0; blue, 0 }  ,draw opacity=0 ][fill={rgb, 255:red, 255; green, 17; blue, 0 }  ,fill opacity=0.33 ][line width=1.5]  (397.56,193.96) -- (617.4,107.25) -- (617.4,336.2) -- (397.56,419.61) -- cycle ;
            \draw  [color={rgb, 255:red, 0; green, 0; blue, 0 }  ,draw opacity=0 ][fill={rgb, 255:red, 12; green, 255; blue, 0 }  ,fill opacity=0.11 ][line width=1.5]  (276.89,105.31) -- (496.73,18.6) -- (496.73,191.14) -- (276.89,274.55) -- cycle ;
            \draw  [color={rgb, 255:red, 0; green, 0; blue, 0 }  ,draw opacity=0 ][fill={rgb, 255:red, 255; green, 17; blue, 0 }  ,fill opacity=0.33 ][line width=1.5]  (276.89,105.31) -- (496.73,18.6) -- (617.4,107.25) -- (397.56,193.96) -- cycle ;
            \draw [line width=1.5]    (276.89,274.55) -- (276.89,108.31) ;
            \draw [shift={(276.89,105.31)}, rotate = 90] [color={rgb, 255:red, 0; green, 0; blue, 0 }  ][line width=1.5]    (14.21,-4.28) .. controls (9.04,-1.82) and (4.3,-0.39) .. (0,0) .. controls (4.3,0.39) and (9.04,1.82) .. (14.21,4.28)   ;
            \draw [line width=1.5]    (222.04,314.85) -- (340.82,460.88) ;
            \draw [shift={(342.71,463.21)}, rotate = 230.88] [color={rgb, 255:red, 0; green, 0; blue, 0 }  ][line width=1.5]    (14.21,-4.28) .. controls (9.04,-1.82) and (4.3,-0.39) .. (0,0) .. controls (4.3,0.39) and (9.04,1.82) .. (14.21,4.28)   ;
            \draw [line width=1.5]    (276.89,274.55) -- (399.37,421.82) ;
            \draw [shift={(401.29,424.12)}, rotate = 230.25] [color={rgb, 255:red, 0; green, 0; blue, 0 }  ][line width=1.5]    (14.21,-4.28) .. controls (9.04,-1.82) and (4.3,-0.39) .. (0,0) .. controls (4.3,0.39) and (9.04,1.82) .. (14.21,4.28)   ;
            \draw [line width=1.5]    (276.89,274.55) -- (493.93,192.21) ;
            \draw [shift={(496.73,191.14)}, rotate = 159.22] [color={rgb, 255:red, 0; green, 0; blue, 0 }  ][line width=1.5]    (14.21,-4.28) .. controls (9.04,-1.82) and (4.3,-0.39) .. (0,0) .. controls (4.3,0.39) and (9.04,1.82) .. (14.21,4.28)   ;
            \draw [line width=1.5]    (222.04,314.85) -- (439.08,232.5) ;
            \draw [shift={(441.88,231.44)}, rotate = 159.22] [color={rgb, 255:red, 0; green, 0; blue, 0 }  ][line width=1.5]    (14.21,-4.28) .. controls (9.04,-1.82) and (4.3,-0.39) .. (0,0) .. controls (4.3,0.39) and (9.04,1.82) .. (14.21,4.28)   ;
            \draw  [color={rgb, 255:red, 0; green, 0; blue, 0 }  ,draw opacity=0 ][fill={rgb, 255:red, 38; green, 249; blue, 50 }  ,fill opacity=0.33 ][line width=0.75]  (276.89,105.31) -- (397.56,193.96) -- (397.56,419.61) -- (276.89,274.55) -- cycle ;
            \draw [line width=1.5]  [dash pattern={on 5.63pt off 4.5pt}]  (276.89,105.31) -- (397.56,193.96) ;
            \draw [line width=1.5]  [dash pattern={on 5.63pt off 4.5pt}]  (276.89,105.31) -- (496.73,18.6) ;
            \draw [line width=1.5]  [dash pattern={on 5.63pt off 4.5pt}]  (617.4,107.25) -- (496.73,18.6) ;
            \draw [line width=1.5]  [dash pattern={on 5.63pt off 4.5pt}]  (617.4,107.25) -- (397.56,193.96) ;
            \draw [line width=1.5]  [dash pattern={on 5.63pt off 4.5pt}]  (397.56,193.96) -- (397.56,419.61) ;
            \draw [line width=1.5]  [dash pattern={on 5.63pt off 4.5pt}]  (617.4,107.25) -- (617.4,336.2) ;
            \draw [line width=1.5]  [dash pattern={on 5.63pt off 4.5pt}]  (617.4,336.2) -- (397.56,419.61) ;
            \draw [line width=1.5]  [dash pattern={on 5.63pt off 4.5pt}]  (496.73,191.14) -- (617.4,336.2) ;
            \draw [line width=1.5]  [dash pattern={on 5.63pt off 4.5pt}]  (135.4,87.6) -- (276.89,105.31) ;
            \draw [line width=1.5]    (222.04,314.85) -- (276.89,274.55) ;
            \draw [line width=1.5]    (134.28,298.73) -- (222.04,314.85) ;
            \draw [line width=1.5]  [dash pattern={on 5.63pt off 4.5pt}]  (342.71,463.21) -- (397.56,419.61) ;
            
            \draw (50,270) node [anchor=north west][inner sep=0.75pt] 
            [xscale=1.6, yscale=1.6] {$(1,1,1)$};
            \draw (35,355.6) node [anchor=north west][inner sep=0.75pt] 
            [xscale=1.6, yscale=1.6] {$(0,0,0)$};
            
            \draw (100,60) node [anchor=north west][inner sep=0.75pt] 
            [xscale=1.6, yscale=1.6] {$(0,0,1)$};
            \draw (245.12,60) node [anchor=north west][inner sep=0.75pt]  
            [xscale=1.6, yscale=1.6]  {$(0,0,1)$};
            
            \draw (394.7,428.13) node [anchor=north west][inner sep=0.75pt]  
            [xscale=1.6, yscale=1.6] {$(1,0,0)$};
            \draw (324.49,468.42) node [anchor=north west][inner sep=0.75pt] 
            [xscale=1.6, yscale=1.6] {$(1,0,0)$};
            
            \draw (497.3,168.83) node [anchor=north west][inner sep=0.75pt]  
            [xscale=1.6, yscale=1.6] {$(0,1,0)$};
            \draw (451.34,227.58) node [anchor=north west][inner sep=0.75pt]  
            [xscale=1.6, yscale=1.6] {$(0,1,0)$};
            
            \draw (150,323.36) node [anchor=north west][inner sep=0.75pt]  
            [xscale=1.6, yscale=1.6] {$(2,2,1)$};
            \draw (190,245) node [anchor=north west][inner sep=0.75pt]  
            [xscale=1.6, yscale=1.6] {$(3,3,2)$};
            \end{tikzpicture}
        }    
        \end{minipage}

        \caption{The tropicalization of the lattice $L$ in \eqref{eq:matrixForExamlpe1}. The polyhedral complex $\Sigma_L$ consists of a line segment, three 2-dimensional polyhedra (in blue) and one 3-dimensional polyhedron (in green). The red regions and dashed lines show where we have truncated the unbounded cells of $\Sigma_L$. The arrows indicate the ray directions.}
        \label{fig:tropL1}
    \end{figure}
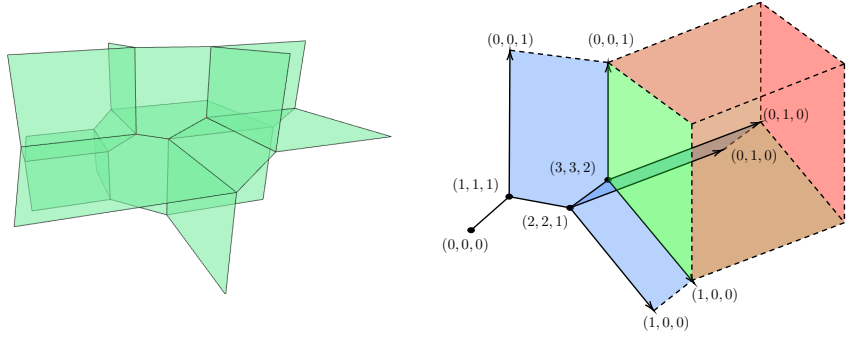
    
    The vectors $u_{I}$ that generate $|\overline{\Sigma}_L|$ over $\overline{\R}_{\geq 0}$ are 
    \[
        \begin{matrix}
                                  & u_{\varnothing} = (0,0,0),   &\\
                    u_{1} = (\infty, 2,1), & u_{2} = (2,\infty,1),  & u_{3} = (1, 1 , \infty), \\
                    u_{12} = (\infty, \infty, 2), &  u_{13} = (\infty, 3 ,\infty) ,          & u_{23} = (3, \infty, \infty)
        \end{matrix}\qquad ,
    \]
    i.e.
        \[
            \trop(L) = |\overline{\Sigma}_L| = \left\{ 
            \begin{bmatrix}
            \min(a,2+c,1+d,3+g)\\  \min(a,2+b,1+d,3+f)\\ \min(a,1+b,1+c,2+e)    
            \end{bmatrix} \colon a,b,c,d,e,f,g \in \overline{\R}_{\geq 0}\right\}.
        \]
    The $2$-skeleton of the polyhedral complex $\mathcal{K}_L$ is the tropical hypersurface defined by the polynomial $\varphi_L$ depicted on the left side ofFigure~\ref{fig:tropL1}. The polyhedral complex $\Sigma_L$ is a subcomplex of $\mathcal{K}_L$ depicted on the right side ofFigure~\ref{fig:tropL1}. 
    \end{example}
    \begin{example}
        Let $L$ be the lattice in $K^3$ spanned over $\Ocal_{K}$ by the rows of the following matrix:
        \begin{equation} \label{eq:matrixForExamlpe2}
            A = \begin{bmatrix}
                    1 & 1   & 1 \\ 
                    0 & t^4 & t^2
            \end{bmatrix}.
        \end{equation} 
        The entropy vector of $L$ is
        \[
        \begin{matrix}
                                   & h_{\varnothing}(L) = 0,   &\\
                     h_{1}(L) = 0, & h_{2}(L) = 0,           & h_{3}(L) = 0, \\
                    h_{12}(L) = 4, & h_{13}(L) = 2,          & h_{23}(L) = 2,\\
                                   & h_{123}(L) = \infty   \ .      &
        \end{matrix}
        \]
        The vectors $u_{I}$ that generate $|\overline{\Sigma}_L|$ over $\overline{\R}_{\geq 0}$ are 
        \[
            \begin{matrix}
                                      & u_{\varnothing} = (0,0,0),   &\\
                        u_{1} = (\infty, 4,2), & u_{2} = (4,\infty,2),  & u_{3} = (2, 2 , \infty)
            \end{matrix}\qquad ,
        \]
        i.e.
        \[
            \trop(L) = |\overline{\Sigma}_L| = \left\{ 
            \begin{bmatrix}
            \min(a,      4+c, 2+d)\\  
            \min(a, 4+b,      2+d)\\ 
            \min(a, 2+b, 2+c     )    
            \end{bmatrix} \colon a,b,c,d \in \overline{\R}_{\geq 0}\right\}.
        \]
        The finite part $\trop(L)\cap\R^3$ is the support of the polyhedral complex $\Sigma_L$ depicted in Figure~\ref{fig:tropL2}.
    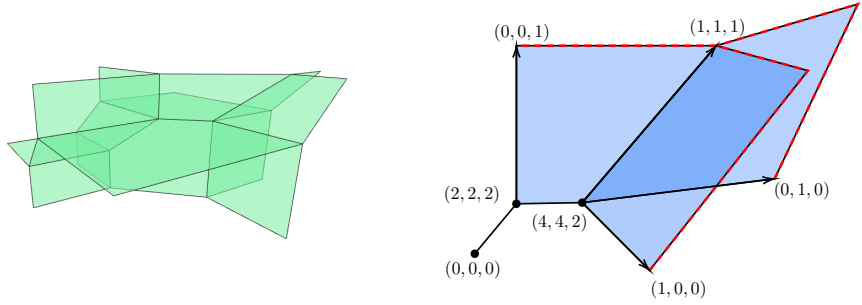
\begin{figure}[ht]
        \centering
        \begin{minipage}{0.45\textwidth}
            \includegraphics[scale=0.27]{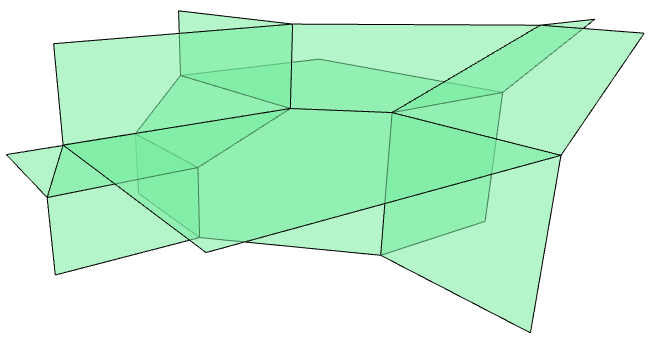}
        \end{minipage}
        \begin{minipage}{0.45\textwidth}
        \scalebox{0.44}{
            \begin{tikzpicture}[x=0.75pt,y=0.75pt,yscale=-1,xscale=1]
            \draw [line width=1.5]    (114,384) -- (164,324) ;
            \draw  [fill={rgb, 255:red, 38; green, 119; blue, 249 }  ,fill opacity=0.33 ][line width=1.5]  (164,134) -- (404,134) -- (243,322.6) -- (164,324) -- cycle ;
            \draw  [fill={rgb, 255:red, 38; green, 119; blue, 249 }  ,fill opacity=0.33 ][line width=1.5]  (404,134) -- (574,84) -- (474,294) -- (243,322.6) -- cycle ;
            \draw  [fill={rgb, 255:red, 38; green, 119; blue, 249 }  ,fill opacity=0.37 ][line width=1.5]  (404,134) -- (514,164) -- (324,404) -- (243,322.6) -- cycle ;
            \draw  [fill={rgb, 255:red, 0; green, 0; blue, 0 }  ,fill opacity=1 ][line width=1.5]  (109.9,384) .. controls (109.9,381.74) and (111.74,379.9) .. (114,379.9) .. controls (116.26,379.9) and (118.1,381.74) .. (118.1,384) .. controls (118.1,386.26) and (116.26,388.1) .. (114,388.1) .. controls (111.74,388.1) and (109.9,386.26) .. (109.9,384) -- cycle ;
            \draw  [fill={rgb, 255:red, 0; green, 0; blue, 0 }  ,fill opacity=1 ][line width=1.5]  (159.9,324) .. controls (159.9,321.74) and (161.74,319.9) .. (164,319.9) .. controls (166.26,319.9) and (168.1,321.74) .. (168.1,324) .. controls (168.1,326.26) and (166.26,328.1) .. (164,328.1) .. controls (161.74,328.1) and (159.9,326.26) .. (159.9,324) -- cycle ;
            \draw  [fill={rgb, 255:red, 0; green, 0; blue, 0 }  ,fill opacity=1 ][line width=1.5]  (238.9,322.6) .. controls (238.9,320.34) and (240.74,318.5) .. (243,318.5) .. controls (245.26,318.5) and (247.1,320.34) .. (247.1,322.6) .. controls (247.1,324.86) and (245.26,326.7) .. (243,326.7) .. controls (240.74,326.7) and (238.9,324.86) .. (238.9,322.6) -- cycle ;
            \draw [color={rgb, 255:red, 255; green, 0; blue, 0 }  ,draw opacity=1 ][line width=2.25]  [dash pattern={on 6.75pt off 4.5pt}]  (404,134) -- (164,134);
            \draw [color={rgb, 255:red, 255; green, 0; blue, 0 }  ,draw opacity=1 ][line width=2.25]  [dash pattern={on 6.75pt off 4.5pt}]  (474,294) -- (574,84) -- (404,134);
            \draw [color={rgb, 255:red, 255; green, 0; blue, 0 }  ,draw opacity=1 ][line width=2.25]  [dash pattern={on 6.75pt off 4.5pt}]  (324,404) -- (514,164) -- (404,134);
            \draw [line width=1.5]    (164,324) -- (164,137);
            
            \draw [shift={(164,134)}, rotate = 90] [color={rgb, 255:red, 0; green, 0; blue, 0 }  ][line width=1.5]    (14.21,-4.28) .. controls (9.04,-1.82) and (4.3,-0.39) .. (0,0) .. controls (4.3,0.39) and (9.04,1.82) .. (14.21,4.28) ;
            \draw [line width=1.5]    (243,322.6) -- (402.05,136.28) ;
            \draw [shift={(404,134)}, rotate = 130.49] [color={rgb, 255:red, 0; green, 0; blue, 0 }  ][line width=1.5]    (14.21,-4.28) .. controls (9.04,-1.82) and (4.3,-0.39) .. (0,0) .. controls (4.3,0.39) and (9.04,1.82) .. (14.21,4.28) ;
            \draw [line width=1.5]    (243,322.6) -- (321.88,401.87);
            \draw [shift={(324,404)}, rotate = 225.14] [color={rgb, 255:red, 0; green, 0; blue, 0 }  ][line width=1.5]    (14.21,-4.28) .. controls (9.04,-1.82) and (4.3,-0.39) .. (0,0) .. controls (4.3,0.39) and (9.04,1.82) .. (14.21,4.28);
            \draw [line width=1.5]    (243,322.6) -- (471.02,294.37);
            \draw [shift={(474,294)}, rotate = 172.94] [color={rgb, 255:red, 0; green, 0; blue, 0 }  ][line width=1.5]    (14.21,-4.28) .. controls (9.04,-1.82) and (4.3,-0.39) .. (0,0) .. controls (4.3,0.39) and (9.04,1.82) .. (14.21,4.28);
            
            \draw (76,390) node [anchor=north west][inner sep=0.75pt] [xscale=1.6, yscale=1.6]  {$(0,0,0)$};
            \draw (75,300) node [anchor=north west][inner sep=0.75pt] [xscale=1.6, yscale=1.6]  {$(2,2,2)$};
            \draw (180,334.4) node [anchor=north west][inner sep=0.75pt] [xscale=1.6, yscale=1.6]   {$(4,4,2)$};
            \draw (135,105) node [anchor=north west][inner sep=0.75pt] [xscale=1.6, yscale=1.6]    {$(0,0,1)$};
            \draw (370,100) node [anchor=north west][inner sep=0.75pt] [xscale=1.6, yscale=1.6]   {$(1,1,1)$};
            \draw (323,414.4) node [anchor=north west][inner sep=0.75pt] [xscale=1.6, yscale=1.6]   {$(1,0,0)$};
            \draw (471,295.4) node [anchor=north west][inner sep=0.75pt] [xscale=1.6, yscale=1.6]   {$(0,1,0)$};
            \end{tikzpicture}
        }
        \end{minipage}
       \caption{The tropicalization of the lattice $L$ in \eqref{eq:matrixForExamlpe2}. The polyhedral complex $\Sigma_L$ consists of a line segment and three 2-dimensional polyhedra. The red regions and dashed lines show where we have truncated the unbounded cells of $\Sigma_L$. The arrows indicate the ray directions.}
        \label{fig:tropL2}
    \end{figure}
    \end{example}

    \noindent A {\tt Julia} \cite{Julia} implementation of the algorithm that computes the polyhedral complex $\Sigma_L$ using {\tt OSCAR} \cite{OSCAR} is available at
    \begin{center}
    {\tt \href{https://www.yelmaazouz.org/TropicalizingLattices.jl}{https://www.yelmaazouz.org/TropicalizingLattices.jl}}.    
    \end{center}

    \subsection{Logarithmic sets of complex polydisks}
        
        We now relate the polyhedral complex $\Sigma_L$ and its support $\trop(L)$ to   logarithmic sets of complex polydisks. 
        Fix positive integers $1 \leq r \leq n$ and let $\mathbb{D} := \{z \in \C \colon |z| \leq 1\}$ denote the closed unit disk in the complex plane. 
        
        Let $A(t)$ be a full-rank $r \times n$ matrix $A(t)$ whose entries are Puiseux series in $t$ such that there exists $\epsilon > 0$ such that $A(t)$ converges uniformly for $|t| \leq \epsilon$. 
        We think of $A$ as an element in $K^{r \times n}$. For series with fractional exponents, we take $t$ positive real when considering the limit.
        We write $L$ for the lattice $L := \Ocal_{K}^r A $ and $\Sigma_L$ is the $r$-dimensional polyhedral complex defined in Definition~\ref{def:SigmaL}. For $0<|t|<1$ we define the map
        \[
        \begin{matrix}
            \Log_t \colon & \C^{n}  & \longrightarrow  &\overline{\R}^{n}, \\
                          & z = (z_1, \dots, z_n)     & \longmapsto   &  \big(\log_{|t|}(|z_1|) , \dots, \log_{|t|}(|z_n|)\big)
        \end{matrix} \quad .
        \]
        For $0<|t|<1$, we denote by $\scrA_t$ the following set in $\R^n$
        \[
            \scrA_t := \Big\{ \Log_t(z) \colon z = y \cdot A(t) \text{ with } y \in \mathbb{D}^r \text{ and } z_i \neq 0 \text{ for all } i=1, \dots, n \Big\}.
        \]
        This is the amoeba of the image of the polydisk $\mathbb{D}^r$ under the linear map $A(t) : \C^r \to \C^n$. We endow $\R^n$ with the euclidean metric. In this setup, the statement of Theorem~\ref{thm:A} can be interpreted as follows: 
        
        \begin{corollary}\label{cor:1}
        The amoeba $\mathscr{A}_t$ converges in Hausdorff distance to  $|\Sigma_L|$, as $t \to 0^+$.
        \end{corollary}
                
        \begin{example}[$n = 2$]
        Consider the family of matrices $A(t) = \begin{bmatrix} 1 & 1 \\ 0 & t\end{bmatrix}$. 
        For $t  \neq 0$ we can write
        \[
            \scrA_t = \left\{ \left( \frac{\log(|\alpha|)}{\log(|t|)} , \ \frac{\log(|\alpha + \beta t|)}{\log(|t|)} \right) : \alpha\in\mathbb{D}\setminus\{0\},\ \beta\in\mathbb{D},\ \alpha+\beta t\neq0 \right\}.
        \]
        Setting $\lambda := -\log(|t|)>0$ and
        $x := \log_{|t|}(|\alpha|)\geq 0$,
        we can describe $\scrA_t$ as
        \[
        \scrA_t =
        \left\{
        (x,y)\in\R^2 :
        \begin{array}{l}
        x\geq 0,\\[2pt]
        y\geq x-\dfrac{1}{\lambda}
        \log(1+e^{\lambda(x-1)}),\\[5pt]
        y\leq x-\dfrac{1}{\lambda}
        \log(1-e^{\lambda(x-1)})
        \quad\text{if }x<1
        \end{array}
        \right\}.
        \]
        This region is depicted inFigure~\ref{fig:amoeba}. 
        We clearly see that, indeed the amoeba $\mathscr{A}_t$ converges to the support of the
        polyhedral complex $\Sigma_L$ when $t^+ \to 0$.

        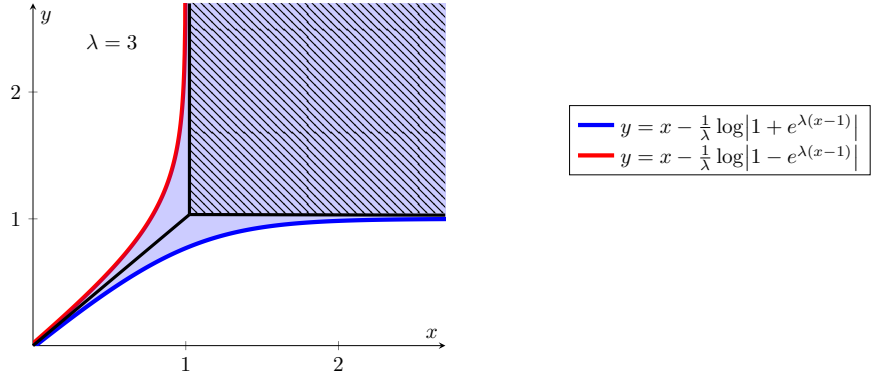
\begin{figure}[H]
            \centering
            \scalebox{0.8}{
                \begin{tikzpicture}
                  \begin{axis}[
                    axis lines=middle,
                    xlabel={$x$},
                    ylabel={$y$},
                    domain=0:2,
                    samples=300,
                    xtick={0,1,2},
                    ytick={0,1,2},
                    ymin=0, ymax=2.7,
                    xmin=0, xmax=2.7,
                    legend style={at={(1.3,0.6)}, anchor=west}
                  ]
                  \addplot[name path = fx, on layer=main,line width = 2, blue, domain=0:6]
                    {x - ln(abs(1 + exp(3*(x - 1)))) / 3};
                  \addlegendentry{$y = x - \frac{1}{\lambda}\log\bigl|1 + e^{\lambda(x - 1)}\bigr|$};
                  \addplot[name path=gx, on layer=main, line width = 2 , red, domain=0:1]
                    {x - ln(abs(1 - exp(3*(x - 1)))) / 3};
                  \addlegendentry{$y = x - \frac{1}{\lambda}\log\bigl|1 - e^{\lambda(x - 1)}\bigr|$};
                  \tikzfillbetween[of = gx and fx, on layer=main]{blue, opacity=0.2};
                  \node[anchor=west] at (axis cs:0.3,2.4)
                    {$\lambda = \pgfmathprintnumber{3}$};
                  \end{axis}
                  \draw[black, line width=1.5] (0, 0)--(2.60, 2.18);
                  \draw[black, line width=1.5] (2.60, 2.18)--(2.60, 5.69);
                  \draw[black, line width=1.5] (2.60, 2.18)--(6.85 , 2.17);
                  \path[pattern color=black, pattern=north west lines] (2.60, 5.7)-- (2.60,2.17) -- (6.85 , 2.2) -- (6.85,5.69) -- cycle;
                \end{tikzpicture}
                }
            \caption{The amoeba $\scrA_t$ (the region colored in blue) converges to the support of the polyhedral complex $\Sigma_L$ (the line segment and shaded orthant) as $t \to 0$ or equivalently as $\lambda \to \infty$ where $\lambda = - \log(|t|)$.}
            \label{fig:amoeba}
        \end{figure}
        
        \end{example}
    
        \section{An algebraic covolume formula}\label{sec:2}

    In this section, we give a purely algebraic interpretation of the entropy
    polynomial.

    \begin{definition}\label{def:covolume}
        Let $M \subset K^r$ be a full-rank $\Ocal_K$-lattice. If the rows
        of $C\in \GL_r(K)$ form an $\Ocal_K$-basis of $M$, set
        \[
            \operatorname{covol}(M):=\val(\det C).
        \]
        Note that this quantity is independent of the chosen basis, since a
        change of basis lies in $\GL_r(\Ocal_K)$ and therefore has
        determinant of valuation zero.
    \end{definition}
    
    The following elementary monotonicity property will be useful in our subsequent arguments.
    
    \begin{lemma}\label{lem:covolumeComparison}
        If $M'\subseteq M$ are full-rank $\Ocal_K$-lattices in $K^r$,
        then we have
        \[
            \operatorname{covol}(M')\geq\operatorname{covol}(M),
        \]
        with equality if and only if $M'=M$.
    \end{lemma}
    \begin{proof}
        Choose basis matrices $C'$ and $C$ for $M'$ and $M$,
        respectively. The inclusion gives $C'=DC$ for some
        $D\in\Mat_{r\times r}(\Ocal_K)$. Hence
        \[
            \operatorname{covol}(M')-\operatorname{covol}(M)
            =\val(\det D)\geq0.
        \]
        Equality holds precisely when $\det D$ is a unit, equivalently
        $D\in\GL_r(\Ocal_K)$, and this is equivalent to $M'=M$.
    \end{proof}

    If $A\in\Mat^\circ_{r\times n}(K)$ and $v\in\R^n$, define the
    \emph{coefficient lattice}
    \begin{equation}\label{eq:coefficientLattice}
        M_v(A):=
        \Big\{y\in\Ocal_K^r:
        \val((yA)_j)\geq v_j\text{ for every }j\in[n]\Big\}.
    \end{equation}
    Note that if $A'=UA$ for $U\in\GL_r(\Ocal_K)$, then $M_v(A')=M_v(A)U^{-1}$. In particular, its covolume does not depend on the choice of the matrix representing $L$. The following formula is the algebraic core for our proofs. The argument is similar to the proof of \cite[Theorem 3.5]{EL22}.

    \begin{proposition}\label{prop:covolumeFormula}
        Let $L=\Ocal_K^r A\subset K^n$, where
        $A\in\Mat^\circ_{r\times n}(K)$. For every $v\in\R^n$, the
        module $M_v(A)$ is a full-rank free $\Ocal_K$-submodule of $\Ocal_K^r$ and we have
        \begin{equation}\label{eq:covolumeFormula}
            \operatorname{covol}(M_v(A))
            =\varphi_L(v)
            =\max_{J\subset[n]}\big(v_J-h_J(L)\big).
        \end{equation}
    \end{proposition}
    \begin{proof}

Set $B := A\operatorname{diag}(t^{-v_1},\dots,t^{-v_n})$. Then we obtain the following
\[
    M_v(A)=\{y\in\Ocal_K^r:yB\in\Ocal_K^n\}.
\]
Using the Smith normal form over the valuation ring $\Ocal_K$, there exist
$V \in \GL_n(\Ocal_K)$ and $U \in \GL_r(\Ocal_K)$ and
$\alpha_1 \leq \dots \leq \alpha_r \in \R$ such that
\[
    UBV=
    \begin{bmatrix}
        t^{\alpha_1}&&&0 & \dots & 0\\
        &\ddots&&\vdots & \ddots & \vdots \\
        &&t^{\alpha_r}&0 & \dots & 0
    \end{bmatrix}.
\]
Since $U$ and $V$ remain invertible over $\Ocal_K$, the change
of variables $z=yU^{-1}$ identifies $M_v(A)$ with the following
lattice
\[
    t^{(-\alpha_1)_+}\Ocal_K\oplus\cdots\oplus
    t^{(-\alpha_r)_+}\Ocal_K,
    \qquad \text{where }
    \quad (a)_+:=\max(a,0).
\]
In particular, $M_v(A)$ is free of rank $r$, and
\begin{equation}\label{eq:covolumeSmith}
    \operatorname{covol}(M_v(A))
    =\sum_{i=1}^r(-\alpha_i)_+.
\end{equation}

We now relate $\alpha_1,\dots,\alpha_r$
to the valuations of the minors of $B$.
For $0\leq k\leq r$, consider
\[
    m_k(B):=
    \min_{\substack{I\subset[r],\,J\subset[n]\\
                    |I|=|J|=k}}
    \val(\det(B_{I\times J})).
\]
We claim that
\[
    m_k(B)=\alpha_1+\cdots+\alpha_k,
\]
where the value for $k=0$ is understood to be zero. Indeed, by the Cauchy--Binet formula, every
$k\times k$ minor of $UBV$ is a sum of products
of a $k\times k$ minor of $B$ with minors of
$U$ and $V$. Since $U$ and $V$ have entries
in $\Ocal_K$, all their minors have nonnegative
valuation. Consequently,
\[
    m_k(UBV)\geq m_k(B).
\]
On the other hand, since $U^{-1}$ and $V^{-1}$
also have entries in $\Ocal_K$, applying the
same argument to gives the reverse inequality. Therefore
\[
    m_k(B) = m_k(UBV).
\]
Now, since $UBV$ is diagonal and $\alpha_1\leq\cdots\leq\alpha_r$,  we obtain as claimed
\[
    m_k(B) = \alpha_1+\cdots+\alpha_k.
\]

By the definition of the entropy vector,
we have
\[
    \min_{\substack{J\subset[n]\\|J|=k}}h_J(B)
    =m_k(B).
\]
Since $h_J(B)=\infty$ whenever $|J|>r$,
it follows that
\[
    \min_{J\subset[n]}h_J(B)
    =\min_{0\leq k\leq r}
      (\alpha_1+\cdots+\alpha_k).
\]
To compute this minimum, recall that
$\alpha_1\leq\cdots\leq\alpha_r$.
The partial sums decrease as long as
$\alpha_k<0$ and increase once
$\alpha_k>0$. Therefore their minimum
is attained when we have added precisely
the negative terms, and therefore we get
\[
    \min_{0\leq k\leq r}
    (\alpha_1+\cdots+\alpha_k)
    =\sum_{\alpha_i<0}\alpha_i.
\]
Combining this with \eqref{eq:covolumeSmith} gives 
$\operatorname{covol}(M_v(A)) = -\min_{J\subset[n]}h_J(B)$. By definition of $B$ we have $h_J(B) = h_J(A) - v_J$. Therefore we obtain the following:
\[
    \operatorname{covol}(M_v(A))
    =-\min_{J\subset[n]}(h_J(L)-v_J)
    =\varphi_L(v),
\]
as claimed.
\end{proof}

    We will also need the following two elementary lemmas for our proofs. 
    The first translates the definition of $\Sigma_L$ into the behavior of the entropy polynomial in the coordinate directions.

    \begin{lemma}\label{lem:Sigma_LOther}
        Let $L$ be a lattice and $h=(h_J)_{J\subset[n]}$ its entropy
        vector. Then
        \[
            |\Sigma_L|
            =
            \big\{v\in\R^n:
            \varphi_L(v+se_j)=\varphi_L(v)+s
            \text{ for all }j\in[n]\text{ and }s>0\big\}.
        \]
        Moreover, if $v \notin |\Sigma_L|$, there are $j\in[n]$ and
        $\varepsilon>0$ such that
        \[
            \varphi_L(v+\varepsilon e_j)=\varphi_L(v).
        \]
    \end{lemma}
    \begin{proof}
        Let
        \[
            \mathcal T_v:=
            \{J\subset[n]:\varphi_L(v)=v_J-h_J\}.
        \]
        Then, by the definition of the labeling of $\Kcal_L$, one has
        $v\in|\Sigma_L|$ if and only if
        $\bigcup_{J\in\mathcal T_v}J=[n]$.
        Suppose first that $v\in|\Sigma_L|$, and fix $j\in[n]$. Some
        $J\in\mathcal T_v$ contains $j$. So, 
        \[
            \varphi_L(v)+s
            =v_J-h_J+s
            \leq\varphi_L(v+se_j)
            \leq\varphi_L(v)+s, \quad \text{for every } s \geq 0,
        \]
        which proves the one inclusion.

        If $v\notin|\Sigma_L|$, choose $j$ belonging to no member of
        $\mathcal T_v$. Since there are only finitely many monomials,
        for all sufficiently small $\varepsilon>0$ no monomial containing
        $j$ catches up with the maximum at $v$. The maximizing monomials
        do not contain $j$, so
        $\varphi_L(v+\varepsilon e_j)=\varphi_L(v)$. 
        This proofs the other inclusion and the second statement.
    \end{proof}

    \begin{lemma}\label{lem:finiteUnionSubmodules}
        Let $R$ be a local ring with infinite residue field, and let $M$
        be a finite free $R$-module. Then $M$ is not a finite union of
        proper submodules. This applies in particular for the case $R = \Ocal_K$
    \end{lemma}
    \begin{proof}
        Let $\mathfrak{m}_R$ be the maximal ideal of $R$, and suppose that
        $N\subsetneq M$ is a submodule. The image
        \[
            (N+\mathfrak{m}_RM) / \mathfrak{m}_RM
            \subset
            M / \mathfrak{m}_RM
        \]
        is a proper linear subspace over the residue field $R/\mathfrak{m}_R$. 
        Indeed, equality would give $M/N=\mathfrak{m}_R(M/N)$, contradicting Nakayama's lemma because
        $M/N$ is a nonzero finitely generated module. A finite-dimensional
        vector space over an infinite field is not a finite union of proper
        linear subspaces. Reducing modulo $\mathfrak{m}_RM$ proves the
        result.
    \end{proof}

    \section{Proof of the main results} \label{sec:3}

    We now prove Theorems~\ref{thm:A} and \ref{thm:B}. As we shall see, the covolume formula in Proposition \ref{prop:covolumeFormula} will be central in our proofs.

    \subsection{Proof of \texorpdfstring{\cref{thm:A}}{Theorem 1.10}}

    \begin{proof}[Proof of Theorem~\ref{thm:A}]
        Write $L=\Ocal_K^r A$ and let $\mathfrak{m}$ be the maximal ideal
        of $\Ocal_K$. We first prove that
        \[
            \trop(L)\cap\R^n=|\Sigma_L|.
        \]
        Let $v=\val(yA)\in\R^n$ for some $y\in\Ocal_K^r$, and suppose
        for a contradiction that $v\notin|\Sigma_L|$. By
       Lemma~\ref{lem:Sigma_LOther}, there are $j\in[n]$ and
        $\varepsilon>0$ such that
        \[
            \varphi_L(v+\varepsilon e_j)=\varphi_L(v).
        \]
        Since $M_{v+\varepsilon e_j}(A)\subseteq M_v(A)$,
        Proposition~\ref{prop:covolumeFormula} gives
        \[
            \operatorname{covol}(M_{v+\varepsilon e_j}(A))
            =\operatorname{covol}(M_v(A)).
        \]
        ByLemma~\ref{lem:covolumeComparison}, we conclude that
        $M_{v+\varepsilon e_j}(A)=M_v(A)$. On the other hand,
        $y\in M_v(A)$ whereas $y\notin M_{v+\varepsilon e_j}(A)$,
        since $\val((yA)_j)=v_j$. This gives a contradiction.

        For the reverse inclusion, fix $v\in|\Sigma_L|$. For every
        $j\in[n]$, consider the submodule
        \[
            N_j:=\{y\in M_v(A):\val((yA)_j)>v_j\}.
        \]
        We claim that $N_j$ is a proper submodule of $M_v(A)$.
        Indeed, suppose that $N_j=M_v(A)$ and choose an
        $\Ocal_K$-basis $b_1,\ldots,b_r$ of $M_v(A)$. Since there are
        finitely many basis vectors, we can find $\varepsilon>0$
        such that
        \[
            \val((b_iA)_j)\geq v_j+\varepsilon
            \quad\text{for all }i\in[r].
        \]
        Every $\Ocal_K$-linear combination satisfies the same
        inequality, and hence $M_{v+\varepsilon e_j}(A)=M_v(A)$.
        But Lemma~\ref{lem:Sigma_LOther} and
        Proposition~\ref{prop:covolumeFormula} give
        \[
            \operatorname{covol}(M_{v+\varepsilon e_j}(A))
            =\varphi_L(v+\varepsilon e_j) =\varphi_L(v)+\varepsilon
              >\operatorname{covol}(M_v(A)),
        \]
        a contradiction. Thus every $N_j$ is proper. Since
        $M_v(A)$ is finite free and $\Ocal_K/\mathfrak{m}=\C$ is
        infinite, Lemma~\ref{lem:finiteUnionSubmodules} gives
        \[
            y\in M_v(A) \Big \backslash \bigcup_{j=1}^nN_j.
        \]
        It follows that $\val((yA)_j)=v_j$ for all $j\in[n]$.
        This proves the equality on the finite part.

        \bigskip
        
        We now pass to the closure in $\overline\R^n$.
        We first show that the vectors $u_J$ appearing in Theorem~\ref{thm:B} are
        valuation vectors of elements of $L$. Fix $J\subsetneq[n]$ with
        $h_J(L)<\infty$ and write $d = |J|$. Since $h_J(L) < \infty$ then if 
        $A$ is a $r \times n$ matrix whose rows are a basis of $L$, there exists $I \subset [r]$ with $|I| = |J|$ and
        \[
            \val(\det(A_{I \times J})) = h_J(L) < \infty.
        \]
        So using elementary row operations, i.e., multiplying $A$ on the left
        by matrices in $\GL_r(\Ocal_K)$, we may assume that $A$ is such that
        \[
            A_{[r]\times J} =
            \begin{bmatrix}C\\0\end{bmatrix},
            \qquad \text{for some} \quad C\in\GL_d(K).
        \]
        The minors of this matrix give
        \[
            h_J(L)=\val(\det C)
        \]
        and, for every $j\notin J$,
        \[
            h_{Jj}(L)
            =h_J(L)+\min_{d<i\leq r}\val(a_{ij}),
        \]
        where the minimum of an empty set is $\infty$.
        Now choose generic coefficients $\beta_{d+1},\dots,\beta_r
        \in\C$ and let
        \[
            x=\sum_{i=d+1}^r\beta_i A_i\in L.
        \]
        Since the residue field is infinite, we may choose the
        coefficients so that the terms of minimum valuation do not
        cancel in any coordinate. Consequently,
        \[
            \val(x_j)=
            \begin{cases}
                \infty &\text{if }j\in J,\\
                h_{Jj}(L)-h_J(L)&\text{if }j\notin J.
            \end{cases}
        \]
        This shows that $\val(x)=u_J$. The same argument works
        for $d=r$. Indeed, the kernel is zero and $u_J=(\infty,\dots,\infty)$
        is the valuation vector of the zero element of $L$.

        Let $C_L$ denote the tropical semimodule generated by
        the vectors $u_J$ over $\overline\R_{\geq0}$. The set
        $\trop(L)$ is stable under nonnegative tropical scalar
        multiplication, since $t^\lambda L\subseteq L$ for
        $\lambda\geq0$. It is also stable under tropical addition:
        for $x,y\in L$, generic constants $\alpha,\beta\in\C$
        satisfy
        \[
            \val(\alpha x+\beta y)=\val(x)\oplus\val(y).
        \]
        In particular $C_L\subseteq\trop(L)$.

        Conversely, let $v\in|\Sigma_L|$. For any proper subset
        $J\subsetneq[n]$ with $h_J(L)<\infty$, set
        \[
            \lambda_J:=\varphi_L(v)+h_J(L)-v_J\geq0.
        \]
        By the definition of the labeling, for each $j\in[n]$
        there is a maximizing monomial indexed by a set containing
        $j$. Therefore,
        \[
            \min_{\substack{J\subsetneq[n]\\h_J<\infty}}
            (\lambda_J+u_{J,j}) 
            =v_j+\varphi_L(v)
              -\max_{\substack{I\subseteq[n]\\j\in I,\ h_I<\infty}}
              (v_I-h_I) = v_j.
        \]
        In the first equality we have used that if $h_I<\infty$
        and $j\in I$, then $h_{I\setminus j}<\infty$.
        This follows from the existence of a nonzero maximal minor
        on the columns indexed by $I$.
        So we deduce that
        \begin{equation}\label{eq:xGenerators}
            v=\bigoplus_{\substack{J\subsetneq[n]\\h_J<\infty}}
              \lambda_J\odot u_J.
        \end{equation}
        In particular, $|\Sigma_L|\subseteq C_L$.

        Since $[0,\infty]$ is compact and tropical addition and
        multiplication by a fixed vector are continuous,
        $C_L$ is closed in $\overline\R^n$. Hence
        \[
            |\overline\Sigma_L|\subseteq C_L\subseteq\trop(L).
        \]
        Finally, let $x\in L$. We can choose $z\in L$ with all
        coordinates nonzero, because $V$ is not contained in any
        coordinate hyperplane and the residue field is infinite.
        For $s\to\infty$, the vectors $x+t^sz\in L$ have all
        coordinates nonzero for sufficiently large $s$, and
        \[
            \val(x+t^sz)\longrightarrow\val(x)
            \quad\text{in }\overline\R^n.
        \]
        By the equality already proved on the finite part, we get
        $\trop(L)\subseteq|\overline\Sigma_L|$. Combining the
        inclusions gives
        \[
            \trop(L)=|\overline\Sigma_L|=C_L.
        \]

        It remains to identify the projectivization. Let
        $\pi:\R^n\to\R^n/\R\bm1_n$ be the quotient map and 
        $V := L\otimes_{\Ocal_K}K$. For any $x\in V\cap(K^\times)^n$,
        write $x=cA$ and choose $s\in\R$ sufficiently large so
        that $t^sc\in\Ocal_K^r$. Then $t^sx\in L$ and
        $\pi(\val(x))=\pi(\val(t^sx))$. Since $K$ is algebraically
        closed and has value group $\R$, the tropical linear
        space $\trop(\PP(V))$ consists of these valuation classes.
        Hence
        \[
            \pi(|\Sigma_L|)=\trop(\PP(V)).
        \]
        Finally, the tropical linear space of the rank-$r$ linear
        subspace $V$ has dimension $r-1$, and $|\Sigma_L|$
        contains the nonnegative $\bm1_n$-ray from each of its
        points. We conclude that $\dim(\Sigma_L)=r$.
    \end{proof}

    \subsection{Proof of \texorpdfstring{\cref{thm:B}}{Theorem 1.12}}

    \begin{proof}[Proof of Theorem~\ref{thm:B}]
        We have shown in the proof of Theorem~\ref{thm:A} that
        $\trop(L)=|\overline\Sigma_L|$ is closed under tropical
        addition and nonnegative tropical scalar multiplication.
        Since $0\in L$, it also contains the tropical zero vector
        $(\infty,\dots,\infty)$, and hence is a tropical
        semimodule over $\overline\R_{\geq0}$.
        Moreover, for every $J\subsetneq[n]$ with $h_J(L)<\infty$,
        we constructed an element $x\in L$ such that
        $\val(x)=u_J$. Thus the semimodule generated by these
        vectors is contained in $\trop(L)$. The reverse inclusion
        follows from \eqref{eq:xGenerators} on the finite part and
        from closedness in $\overline\R^n$. Therefore
        \[
            |\overline\Sigma_L|
            =\bigoplus_{\substack{J\subsetneq[n]\\h_J(L)<\infty}}
              \overline\R_{\geq0}\odot u_J,
        \]
        and this concludes the proof.
    \end{proof}

    \section{Beyond the realizable case}
    
    In tropical geometry, tropical linear spaces arise from tropicalizing linear subspaces over a non-archimedean field or more generally from the broader combinatorial framework of valuated matroids. In our setting, instead of tropicalizing non-archimedean lattices, we consider valuated bimatroids in the sense of \cite{Murota} and define their entropy vectors in parallel with the realizable case. Following the notation of \cite[Definition 2.1]{GRSU}, we adopt the following definition of valuated bimatroids.

    For any pair of integers $n,r \geq 1$ we write $\binom{[r]}{\ast}\times\binom{[n]}{\ast}$ for the set
    \[
    \binom{[r]}{\ast}\times\binom{[n]}{\ast} := \{ (I,J) \colon I \subset [r], J \subset [n] \text{ and } |I| =|J| \}.
    \]
    We use the notation $Ii := I\cup\{i\}$ and $I\setminus i := I\setminus\{i\}$.

    \begin{definition}
    
    Let $n,r \geq 1$ be integers. 
    A valuated bimatroid ${\mathrm{A}}$ with rows indexed by $[r]$ and 
    columns indexed by $[n]$ is a map
        \[
            \nu_{\mathrm{A}} \colon \binom{[r]}{\ast} \times \binom{[n]}{\ast} \longrightarrow \overline{\R}, \quad (I, J) \mapsto \nu_{\mathrm{A}}(I,J)
        \]
        satisfying the following list of axioms:
        \begin{enumerate}[wide = 0pt]
            \item \label{item:(1)} $\nu_{\mathrm{A}}(\varnothing, \varnothing) = 0$,

            \item for any $(I,J), (I', J') \in \binom{[r]}{\ast}\times \binom{[n]}{\ast}$ the following holds

            \begin{enumerate}[wide=20pt]
                \item \label{item:(2a)} for any $i' \in I' \setminus I$ at least one of the following statements holds:
                
                \begin{itemize}[wide = 40pt]
                    \item there exists $i \in I \setminus I'$ such that
                    \[
                        \nu_{\mathrm{A}} (I,J)  +  \nu_{\mathrm{A}}(I', J') \geq \nu_{\mathrm{A}} \Big( Ii' \setminus i , J\Big) +  \nu_{\mathrm{A}} \Big(I'i \setminus i', J' \Big).
                    \]

                    \item there exists $j' \in J' \setminus J$ such that
                   \[
                        \nu_{\mathrm{A}} (I,J)  +  \nu_{\mathrm{A}}(I', J') \geq \nu_{\mathrm{A}} \Big(Ii', Jj'\Big) +  \nu_{\mathrm{A}} \Big(I' \setminus i' , J' \setminus j' \Big).
                    \]

                \end{itemize}

                \item \label{item:(2b)}  for any $j \in J \setminus J'$ at least one of the following two statements holds:
                
                \begin{itemize}[wide = 40pt]
                    \item  there exists $i \in I\setminus I'$ such that
                    \[
                        \nu_{\mathrm{A}} (I,J)  +  \nu_{\mathrm{A}}(I', J') \geq \nu_{\mathrm{A}} \Big(I \setminus i , J \setminus j \Big) +  \nu_{\mathrm{A}} \Big(I'i, J'j\Big).
                    \]
                    \item  there exists $j' \in J'\setminus J$ such that
                    \[
                        \nu_{\mathrm{A}} (I,J)  +  \nu_{\mathrm{A}}(I', J') \geq \nu_{\mathrm{A}} \Big(I, J j' \setminus  j \Big) +  \nu_{\mathrm{A}} \Big(I', J'j \setminus j'\Big).
                    \]
                    
                \end{itemize}
            \end{enumerate}

        \end{enumerate}
    \end{definition}

    \begin{definition}
        Given a valuated bimatroid ${\mathrm{A}}$ on $\binom{[r]}{\ast} \times \binom{[n]}{\ast}$, we define its entropy vector $h({\mathrm{A}})$ as follows:
        \[
         h_{J}(\mathrm{A}) := \min_{I \in \binom{[r]}{|J|}} \nu_{\mathrm{A}}(I,J) \quad \text{for any } J \subset [n].
        \]
        Similar to Definition~\ref{def:SigmaL}, we define the polyhedral complex $\Sigma_{\mathrm{A}}$ through the \emph{entropy polynomial} 
        \[
        \varphi_{\mathrm{A}}(x) := 
        \max_{J \subset [n]}(x_J - h_J(\mathrm{A})).
        \]
        We write $|\overline{\Sigma}_{\mathrm A}|$ for the closure
        of $|\Sigma_{\mathrm A}|$ in $\overline{\R}^{n}$.
    \end{definition}
    
\begin{definition}[cf.~\cite{Murota}]
    A function $h:2^{[n]}\to\overline{\mathbb{R}}$ is called
    \emph{$M^\natural$-convex} if for any subsets
    $I,J\subseteq[n]$ and any $i\in I\setminus J$,
    at least one of the following holds:
    \begin{enumerate}[label=(E\arabic*),ref=(E\arabic*)]
        \item \label{item:MnExchange1}
        $h_I + h_J \geq h_{I\setminus i} + h_{Ji}$.

        \item \label{item:MnExchange2}
        There exists $j\in J\setminus I$ such that
        $h_I + h_J \geq h_{Ij\setminus i} + h_{Ji\setminus j}$.
    \end{enumerate}
\end{definition}

We note that $M^\natural$-convexity is stronger property than supermodularity, as the following proposition shows. This is a consequence of \cite[Theorem 6.19]{Murota2}, but for completeness we include a proof using our notation.

\begin{proposition}[Theorem 6.19 in \cite{Murota2}] \label{prop:Supermod}
    Any $M^\natural$-convex function $h:2^{[n]}\to\overline{\mathbb{R}}$ is supermodular.
\end{proposition}
\begin{proof}
Let $h:2^{[n]}\to\overline{\mathbb{R}}$ be an $M^\natural$-convex function.
Let $S\subseteq[n]$ and let $i,j\notin S$ be distinct. 
Applying the $M^\natural$-exchange axiom to $Sij$ and $S$,
the second alternative is impossible since $S\setminus Sij=\varnothing$. Therefore
\[
    h_{Sij}+h_S\geq h_{Si}+h_{Sj}.
\]
In particular, for any $S\subseteq T$ and
$i\notin T$, whenever $h_{Ti}<\infty$,
we obtain
\[
    h_{Si}-h_S\leq h_{Ti}-h_T
\]
by adding the elements of $T\setminus S$
one at a time.

Now fix $I,J\subseteq[n]$. If
$h_{I\cup J}=\infty$, there is nothing
to prove. Otherwise, all the subsets of
$I\cup J$ have finite entropy by (C1).
Writing $J\setminus I=\{j_1,\dots,j_m\}$
and applying the preceding inequality
successively gives
\[
    h_J-h_{I\cap J}
    \leq h_{I\cup J}-h_I.
\]
Thus
\[
    h_I+h_J
    \leq h_{I\cap J}+h_{I\cup J},
\]
which proves supermodularity.
\end{proof}

    \begin{theorem}
        For any valuated bimatroid $\mathrm{A}$ on $\binom{[r]}{\ast} \times \binom{[n]}{\ast}$ the following hold:

        \begin{enumerate}
            \item 
            The entropy vector $h({\rm A}) = (h_{J}(\mathrm{A}))_{J \subset [n]}$ is $M^\natural$-convex. In particular it is supermodular i.e.
            \[
                h_I(\mathrm{A}) + h_{J}(\mathrm{A}) \leq  h_{I \cap J}(\mathrm{A}) + h_{I \cup J}(\mathrm{A}) \quad \text{for any } I,J \subset [n].
            \]
            
            \item Suppose that $\mathrm{A}$ is loopless i.e. that $h_j(\mathrm{A}) < +\infty $ for all $j \in [n]$. Then the support $|\overline{\Sigma}_{\mathrm{A}}|$ of the polyhedral complex $\overline{\Sigma}_{\mathrm{A}}$ is a tropical semimodule over $\overline{\R}_{\geq 0}$ generated by the vectors $\{ u_{J} \colon J \subsetneq [n], \ h_{J} \neq \infty \}$ in $\overline{\R}^n$ defined as follows 
        \[
            u_{J,j} := \begin{cases}
                        h_{Jj}(\mathrm{A}) - h_{J}(\mathrm{A}) & \text{if } j \not \in J \\ 
                             +  \infty & \text{otherwise} 
                      \end{cases}, \quad \text{where } Jj := J \cup \{j\}.
        \]
   
        \end{enumerate}
        
    \end{theorem}

\begin{proof}
    To simplify notation, we write $h_J=h_J(\mathrm A)$, and let
    \[
        \mathcal{I}:=\{J\subseteq[n]:h_J<\infty\}.
    \]

    We first prove that $h(\mathrm A)$ is
    $M^\natural$-convex. Let $J,J'\in\mathcal{I}$
    and fix $j\in J\setminus J'$.
    Choose $I,I'\subseteq[r]$ such that
    \[
        \nu_{\mathrm A}(I,J)=h_J,
        \qquad
        \nu_{\mathrm A}(I',J')=h_{J'}.
    \]
    Applying axiom \ref{item:(2b)} to
    $(I,J)$ and $(I',J')$, we obtain one
    of the following two inequalities.

    \begin{itemize}
        \item In the first case, there exists
    $i\in I\setminus I'$ such that
    \[
        h_J+h_{J'}
        \geq
        \nu_{\mathrm A}(I\setminus i,J\setminus j)
        +\nu_{\mathrm A}(I'i,J'j) \geq h_{J\setminus j}+h_{J'j}.
    \]
    \item In the second case, there exists
    $j'\in J'\setminus J$ such that
    \[
        h_J+h_{J'}
        \geq
        \nu_{\mathrm A}(I,Jj'\setminus j)
        +\nu_{\mathrm A}(I',J'j\setminus j')
        \geq
        h_{Jj'\setminus j}+h_{J'j\setminus j'}.
    \]
    \end{itemize}
    These are precisely the two alternatives (\ref{item:MnExchange1}-\ref{item:MnExchange2}) in the definition of $M^\natural$-convexity. 
    If either $h_J$ or $h_{J'}$ is infinite,
    the exchange inequality is automatic.
    Thus $h(\mathrm A)$ is $M^\natural$-convex. 
    By Proposition \ref{prop:Supermod}, we deduce that
    $h(\mathrm A)$ is also supermodular. This proves the first assertion.

    \medskip

    To prove the second assertion, we now assume that $\mathrm A$ is loopless. 
    We will use the following two consequences of the exchange property (\ref{item:MnExchange1} and \ref{item:MnExchange2}).

    \begin{enumerate}
    \item  
    First, $\mathcal{I}$ is closed under taking
    subsets. Indeed, if $J\in \mathcal{I}$
    and $j\in J$, applying the
    $M^\natural$-exchange property \ref{item:MnExchange1} to
    $J$ and $\varnothing$ gives
    \[
        h_J+h_\varnothing
        \geq h_{J\setminus j}+h_j,
    \]
    since the exchange alternative \ref{item:MnExchange2} is
    impossible. As $h_\varnothing=0$,
    this implies $h_{J\setminus j}<\infty$.

    \item Second, for any $J,J'\in \mathcal{I}$
    and $j\in J\setminus J'$, we have either (alternative \ref{item:MnExchange1})
    \begin{equation}\label{eq:projectedExchange1}
        h_J+h_{J'}
        \geq h_{J\setminus j}+h_{J'j},
    \end{equation}
    or (alternative \ref{item:MnExchange2}) there exists $j'\in J'\setminus J$
    such that
    \begin{equation}\label{eq:projectedExchange2}
        h_J+h_{J'}
        \geq
        h_{Jj'\setminus j}
        +h_{J'j\setminus j'}.
    \end{equation}
    \end{enumerate}
    For $x\in\R^n$, set
    \[
        \ell_x(J):=x_J-h_J,
        \qquad
        \mathcal B(x):=
        \{J\subseteq[n]:
        \ell_x(J)=\varphi_{\mathrm A}(x)\}.
    \]
    By the definition of $\Sigma_{\mathrm A}$,
    \begin{equation}\label{eq:supportdef}
        |\Sigma_{\mathrm A}|
        =
        \left\{
        x\in\R^n:
        \bigcup_{J\in\mathcal B(x)}J=[n]
        \right\}.
    \end{equation}

    We first show that $|\Sigma_{\mathrm A}|$
    is closed under tropical addition and
    nonnegative tropical scalar multiplication. To do so, we shall need the following lemma.

    \begin{lemma}
    We claim that if $y\leq x$ in $\R^n$
    and $y_j=x_j$ for some
    $j\in J\in\mathcal B(x)$, then
    there exists $J'\in\mathcal B(y)$
    containing $j$.
    \end{lemma}
    \begin{proof}[Proof of Lemma]
        To prove this, choose
    $J'\in\mathcal B(y)$.
    If $j\in J'$, there is nothing to show.
    Otherwise, apply the exchange property \ref{item:MnExchange1}--\ref{item:MnExchange2}
    to $J,J'$ and $j$. 
    
    \begin{itemize}
        \item In the first case \ref{item:MnExchange1}, we have
    \(
        h_J+h_{J'}
        \geq h_{J\setminus j}+h_{J'j},
    \)
    and hence
    \[
        \ell_x(J)+\ell_x(J')
        \leq
        \ell_x(J\setminus j)+\ell_x(J'j).
    \]
    Since $J\in\mathcal B(x)$, we obtain
    \(
        \ell_x(J')\leq\ell_x(J'j).
    \)
    Since $x_j=y_j$, this gives
    \(
        \ell_y(J')\leq\ell_y(J'j),
    \)
    and hence $J'j\in\mathcal B(y)$.

    \item In the second case \ref{item:MnExchange2}, there exists
    $k\in J'\setminus J$ such that
    \[
        h_J+h_{J'}
        \geq
        h_{Jk\setminus j}+h_{J'j\setminus k}.
    \]
    Arguing as above, we obtain
    \[
        \ell_x(J')
        \leq\ell_x(J'j\setminus k).
    \]
    Since $x_j=y_j$ and $y_k\leq x_k$,
    \[
        \ell_y(J'j\setminus k)-\ell_y(J')
        =
        \ell_x(J'j\setminus k)-\ell_x(J') +(x_k-y_k) \geq 0.
    \]
    Therefore $J'j\setminus k\in\mathcal B(y)$.
    \end{itemize}
    
    This proves the lemma.
    \end{proof}

    Now let $x,y\in|\Sigma_{\mathrm A}|$
    and set $z = x\oplus y$.
    Fix $j\in[n]$. Without loss of generality we may assume that $x_j = z_j$.
    We have $z\leq x$, and using \eqref{eq:supportdef} we deduce that $j$ belongs to some $J \in \mathcal B(x)$.
    By the preceding claim, $j$ belongs to some member of $\mathcal B(z)$.
    Since this holds for every $j \in [n]$,
    \eqref{eq:supportdef} gives
    \[
        z = x \oplus y \in |\Sigma_{\mathrm A}|.
    \]

    Next, let $x\in|\Sigma_{\mathrm A}|$,
    let $\lambda\geq0$, and set
    $z=x+\lambda\bm1_n$.
    Fix $j\in[n]$ and choose $J\in\mathcal B(x)$ containing $j$
    and $J'\in\mathcal B(z)$.
    If $j\notin J'$, we apply \ref{item:MnExchange1}--\ref{item:MnExchange2} to $J,J'$ and~$j$.
    \begin{itemize}
        \item In the first case (alternative \ref{item:MnExchange1}), the same argument
    as above gives
    \[
        \ell_x(J')\leq\ell_x(J'j).
    \]
    Consequently,
    \[
    \begin{aligned}
        \ell_z(J'j)-\ell_z(J')
        =
        \ell_x(J'j)-\ell_x(J')+\lambda
        \geq0.
    \end{aligned}
    \]
    Thus $J'j\in\mathcal B(z)$.

    \item In the second case (alternative \ref{item:MnExchange2}), there exists
    $k\in J'\setminus J$ such that
    \[
        \ell_x(J')
        \leq\ell_x(J'j\setminus k).
    \]
    Since $J'$ and $J'j\setminus k$
    have the same cardinality,
    \[
        \ell_z(J'j\setminus k)-\ell_z(J')
        =
        \ell_x(J'j\setminus k)-\ell_x(J')
        \geq0.
    \]
    Hence $J'j\setminus k\in\mathcal B(z)$.
    \end{itemize}
    We conclude that
    \[
        \lambda\odot x
        =x+\lambda\bm1_n
        \in|\Sigma_{\mathrm A}|.
    \]
    Thus $|\Sigma_{\mathrm A}|$ is closed
    under tropical addition and tropical
    scalar multiplication by elements of
    $\R_{\geq0}$. 
    
    We now show that the vectors $u_J$ belong to $|\overline\Sigma_{\mathrm A}|$. First, consider the case $J = \varnothing$.
    Since $h_\varnothing = 0$,
    supermodularity gives
    \[
        h_I\geq\sum_{i\in I}h_i
        \qquad\text{whenever }h_I<\infty.
    \]
    At $u_\varnothing=(h_1,\dots,h_n)$,
    we therefore have
    \[
        \ell_{u_\varnothing}(I)\leq0
        \qquad\text{for all }I\subseteq[n].
    \]
    Moreover, the empty set and every
    singleton are maximizing, since
    $\mathrm A$ is loopless.
    Thus
    \[
        u_\varnothing\in|\Sigma_{\mathrm A}|.
    \]
    
    Note that $ \mathcal{I}$
    is the collection of independent sets
    of a matroid $\mathscr{M}$ on $[n]$. Indeed,

    \begin{enumerate}
        \item   $\varnothing \in \mathcal{I}$
        
        \item   we have already shown that $\mathcal{I}$ is closed under
                taking subsets.
        
        \item   let $I,I'\in \mathcal{I}$ with
                $|I| > |I'|$ and choose
                $i\in I\setminus I'$.
                The exchange property \ref{item:MnExchange1}--\ref{item:MnExchange2} implies that
                either $I'i\in \mathcal{I}$, or
                there exists $j\in I'\setminus I$
                such that $I'i\setminus j\in \mathcal{I}$.
                In the second case (alternative \ref{item:MnExchange2}), replacing $I'$
                by $I'i\setminus j$ decreases
                $|I'\setminus I|$ while preserving
                its cardinality.
                Repeating this procedure eventually
                gives an independent set contained
                in $I$ to which an element of $I$
                can be added. This proves the augmentation axiom.
    \end{enumerate}

    Now fix $J\in \mathcal{I}$ with $0 < d := |J| < n$. We claim that 
    \[
        u_J \in |\Sigma_{\rm A}|.
    \]
    To prove this, we shall construct a sequence of points $x^{(T)}$ such that $\lim\limits_{T \to \infty} x^{(T)} = u_J$.
    
    \begin{enumerate}
    \item [\underline{\bf Step 1})] 
    First let $C$ be the closure of $J$ in the matroid $\mathscr{M}$, namely
    \[
        C:= J \cup \{i\in[n]\setminus J:h_{Ji}=\infty\}
    \]
    of equivalently for any $j \in [n]$ we have
    \[
        u_{J,j} = \infty \quad \text{ if and only if } \quad j \in C.
    \]
    Then $J$ is a basis of $C$.  For $i\in J$, set $p_i=0$ and for $i\in C\setminus J$ define
    \[
        p_i:=
        \min_{\substack{j\in J\\
                    Ji\setminus j\in \mathcal{I}}}
        (h_{Ji\setminus j}-h_J).
    \]
    Note that the minimum here is finite by the basis
    exchange property. We claim that
    \begin{equation}\label{eq:bimatroidBasisWeights}
        p_B-h_B\leq-h_J, \quad \text{for every basis } B \text{ of } C.    \end{equation}
    Indeed, if $B=J$, equality holds.
    Otherwise, choose $i\in B\setminus J$
    and apply the exchange property \ref{item:MnExchange1}--\ref{item:MnExchange2}
    to $B,J$ and $i$.
    Since $Ji$ is dependent, the first
    alternative \ref{item:MnExchange1} is impossible.
    Consequently, by alternative \ref{item:MnExchange2}, there exists
    $j\in J\setminus B$ such that
    \[
        h_B+h_J
        \geq
        h_{Bj\setminus i}+h_{Ji\setminus j}.
    \]
    By the definition of $p_i$,
    \[
        h_{Ji\setminus j}-h_J\geq p_i,
    \]
    and hence
    \[
        h_B-h_{Bj\setminus i}\geq p_i.
    \]
    Since $p_j=0$, we obtain
    \[
        p_B-h_B
        \leq
        p_{Bj\setminus i}-h_{Bj\setminus i}.
    \]
    Repeating this argument decreases
    $|B\setminus J|$ until we reach $J$,
    proving \eqref{eq:bimatroidBasisWeights}. 
    Moreover, for every
    $i\in C\setminus J$, there exists
    $j\in J$ attaining the minimum
    in the definition of $p_i$.
    For $B=Ji\setminus j$, we then have
    \[
        p_B-h_B=-h_J.
    \]

    \item [\underline{\bf Step 2})] 
    We next show that for every basis
    $B$ of $C$ and every $k\notin C$,
    \begin{equation}\label{eq:bimatroidMarginalInvariant}
        h_{Bk}-h_B=h_{Jk}-h_J.
    \end{equation}
    It is enough to prove that the quantity
    $h_{Bk}-h_B$ is unchanged under
    a single basis exchange in $C$,
    since the graph of bases of $C$
    is connected.
    Let $B'=Bi\setminus j$ be another
    basis of $C$, where
    $j\in B$ and $i\in C\setminus B$.
    Applying the exchange property \ref{item:MnExchange1}--\ref{item:MnExchange2}
    to the pair $Bk,B'$ and the element $k$, 
    the exchange alternative \ref{item:MnExchange2} is impossible
    because $Bi$ is dependent.
    Therefore, by alternative \ref{item:MnExchange1},
    \[
        h_{Bk}+h_{B'}
        \geq h_B+h_{B'k}.
    \]
    Interchanging $B$ and $B'$ gives
    the reverse inequality, and hence
    \[
        h_{Bk}-h_B
        =
        h_{B'k}-h_{B'}.
    \]
    This proves
    \eqref{eq:bimatroidMarginalInvariant}.

    \item [\underline{\bf Step 1})] 
    We now construct the sequence of vectors $x^{(T)}$. 
    For $k\notin C$, we set
    \[
        w_k:=h_{Jk}-h_J,
    \]
    and for $T>0$, we define $x^{(T)}\in\R^n$ as follows
    \[
        x^{(T)}_i:=
        \begin{cases}
            T+p_i & \text{if }i\in C,\\
            w_i & \text{if }i\notin C.
        \end{cases}
    \]
    We claim that
    $x^{(T)}\in|\Sigma_{\mathrm A}|$
    for all sufficiently large $T$. Indeed, since $J$ is a basis of $C$,
    every $I\in \mathcal{I}$ satisfies
    $|I\cap C|\leq d$.
    
\begin{itemize}
    \item If $|I\cap C|<d$, then
    \[
        \ell_{x^{(T)}}(I) = \sum_{i \in I} x^{(T)}_i - h_I
        <
        dT-h_J
        =
        \ell_{x^{(T)}}(J)
    \]
    for all sufficiently large $T$.

    \item   
    Now suppose $|I\cap C|=d$.
    Write $I=B\cup U$, where
    $B$ is a basis of $C$ and
    $U\subseteq[n]\setminus C$.
    By supermodularity,
    \[
        h_{B\cup U}
        \geq
        h_B+
        \sum_{k\in U}(h_{Bk}-h_B).
    \]
    Using
    \cref{eq:bimatroidBasisWeights,eq:bimatroidMarginalInvariant},
    we obtain
    \[
        \ell_{x^{(T)}}(I)
        =
        dT+p_B+\sum_{k\in U}w_k-h_I
        \leq
        dT + p_B - h_B
        \leq dT-h_J.
    \]
    
    \end{itemize}
    Thus $J$ is a maximizing monomial of $\varphi_{\mathrm A}$ at $x^{(T)}$. Furthermore, for every $i\in C\setminus J$, there exists $j\in J$ such that
    \[
        p_i=h_{Ji\setminus j}-h_J.
    \]
    Consequently,
    \[
        \ell_{x^{(T)}}(Ji\setminus j)
        =
        dT-h_J,
    \]
    so $Ji\setminus j$ is also maximizing.
    Similarly, for every $k\notin C$,
    \[
        \ell_{x^{(T)}}(Jk)
        =
        dT+w_k-h_{Jk}
        =
        dT-h_J,
    \]
    and hence $Jk$ is a maximizing set. The maximizing subsets therefore cover $[n]$. So by \eqref{eq:supportdef}, we deduce that
    \[
        x^{(T)}\in|\Sigma_{\mathrm A}| \quad \text{ for all sufficiently large } T.
    \]
    Finally, by the definition of $C$,
    \[
        \lim_{T\to\infty}x^{(T)}_i
        =
        \begin{cases}
            \infty & \text{if }i\in C,\\
            h_{Ji}-h_J & \text{if }i\notin C.
        \end{cases}
    \]
    This is precisely the vector $u_J$.
    \end{enumerate}
Therefore, as claimed, we have:
\[
    u_J \in|\overline\Sigma_{\mathrm A}|.
\]

Since $u_\varnothing\in|\Sigma_{\mathrm A}|$
and $|\Sigma_{\mathrm A}|$ is closed
under nonnegative tropical scalar
multiplication, the all-infinite vector
also belongs to its closure.
Tropical addition and scalar
multiplication extend continuously to
$\overline\R^n$ and
$\overline\R_{\geq0}$.
Hence $|\overline\Sigma_{\mathrm A}|$
is a tropical semimodule over
$\overline\R_{\geq0}$.
As every $u_J$ belongs to this
semimodule, we then have
\[
    \bigoplus_{\substack{J\subsetneq[n]\\
                        h_J<\infty}}
    \overline\R_{\geq0}\odot u_J
    \subseteq
    |\overline\Sigma_{\mathrm A}|.
\]

It remains to prove the reverse inclusion. To that end, 
fix $x\in|\Sigma_{\mathrm A}|$.
For every $J\subsetneq[n]$
with $h_J<\infty$, set
\[
    \lambda_J:=
    \varphi_{\mathrm A}(x)+h_J-x_J.
\]
By the definition of
$\varphi_{\mathrm A}$, we have
$\lambda_J\geq0$. For every $j\in[n]$ we have
\[
    \min_{\substack{J\subsetneq[n]\\
                    h_J<\infty}}(\lambda_J+u_{J,j})
    = \min_{\substack{J\subsetneq[n]\\ j\notin J,\ h_J<\infty}}
    \bigl(
    \varphi_{\mathrm A}(x) +h_{Jj}-x_J \bigr)
    =
    x_j+\varphi_{\mathrm A}(x) - \max_{\substack{I\subseteq[n]\\
                    j\in I,\ h_I<\infty}}
    (x_I-h_I).
\]
Here we used the fact that every
$I\in \mathcal{I}$ containing $j$
can be written as
$I=(I\setminus j)j$, with
$I\setminus j\in \mathcal{I}$.
Since $x\in|\Sigma_{\mathrm A}|$,
there exists a maximizing subset
containing $j$. Therefore
\[
    \max_{\substack{I\subseteq[n]\\
                    j\in I,\ h_I<\infty}}
    (x_I - h_I)
    =
    \varphi_{\mathrm A}(x),
\]
and consequently
\[
    \min_{\substack{J\subsetneq[n]\\
                    h_J<\infty}}
    (\lambda_J+u_{J,j})
    =
    x_j.
\]
Since this holds for every $j\in[n]$,
\[
    x=
    \bigoplus_{\substack{J\subsetneq[n]\\
                        h_J<\infty}}
    \lambda_J\odot u_J.
\]
Thus every point of $|\Sigma_{\mathrm A}|$
belongs to the tropical semimodule
generated by the vectors $u_J$.
This semimodule is closed in
$\overline\R^n$, since it is the
continuous image of a finite product
of copies of $[0,\infty]$.
Taking the closure gives
\[
    |\overline\Sigma_{\mathrm A}|
    \subseteq
    \bigoplus_{\substack{J\subsetneq[n]\\
                        h_J<\infty}}
    \overline\R_{\geq0}\odot u_J.
\]
Combining the two inclusions,
we conclude that
\[
    |\overline\Sigma_{\mathrm A}|
    =
    \bigoplus_{\substack{J\subsetneq[n]\\
                        h_J<\infty}}
    \overline\R_{\geq0}\odot u_J,
\]
as claimed.
\end{proof}

    \section*{Acknowledgments}
    This work benefited from valuable discussions with Shelby Cox, Bernd
    Sturmfels, and Cynthia Vinzant, to whom I am very grateful. I also thank
    Weihong Xu for proofreading an earlier version of the paper, and Yue Ren
    for assistance with computations in {\tt OSCAR}.

    \section*{AI Disclosure}
    AI-assisted tools were used exclusively for editorial and typesetting
    purposes, including identifying typographical errors in the manuscript.
    They were not used in the development of the mathematical arguments, 
    proofs, or results.

\end{document}